\documentclass[12pt]{amsart}
\usepackage{amsmath,amsthm}
\usepackage{amssymb}
\usepackage{euscript}
\DeclareMathAlphabet{\EuRm}{U}{eur}{m}{n}
\SetMathAlphabet{\EuRm}{bold}{U}{eur}{b}{n}
\begin{document}
%
%
\swapnumbers
\newtheorem{thm}{Theorem}[section]
\newtheorem*{tha}{Theorem A}
\newtheorem*{thb}{Theorem B}
\newtheorem{lemma}[thm]{Lemma}
\newtheorem{prop}[thm]{Proposition}
\newtheorem{cor}[thm]{Corollary}
\theoremstyle{definition}
\newtheorem{defn}[thm]{Definition}
\newtheorem{example}[thm]{Example}
\newtheorem{summary}[thm]{Summary}
\newtheorem{fact}[thm]{Fact}
\theoremstyle{remark}
\newtheorem{notation}[thm]{Notation}
\newtheorem{remark}[thm]{Remark}

\newtheorem{note}[thm]{Note}
\newtheorem{aside}[thm]{Aside}
\numberwithin{equation}{section}
%
%
\def\sect{\setcounter{thm}{0}\section}
%
%
\newcommand{\hsp}{\hspace{10 mm}}
\newcommand{\hs}{\hspace{5 mm}}
\newcommand{\hsm}{\hspace{2 mm}}
\newcommand{\vs}{\vspace{7 mm}}
\newcommand{\vsm}{\vspace{2 mm}}
%
%
\newcommand{\xra}[1]{\xrightarrow{#1}}
\newcommand{\xla}[1]{\xleftarrow{#1}}
\newcommand{\hra}{\hookrightarrow}
\newcommand{\rest}[1]{\lvert_{#1}}
\newcommand{\lra}[1]{\langle{#1}\rangle}
\newcommand{\lrv}[1]{\lvert{#1}\rvert}
\newcommand{\DEF}{:=}
\newcommand{\EQUIV}{\Leftrightarrow}
\newcommand{\IMPLY}{\Rightarrow}
\newcommand{\epic}{\to\hspace{-5 mm}\to}
\newcommand{\xepic}[1]{\xrightarrow{#1}\hspace{-5 mm}\to}
\newcommand{\vars}{\varsigma}
%
%
\newcommand{\Coker}{\operatorname{Coker}}
\newcommand{\colim}{\operatorname{colim}}
\newcommand{\cosk}[1]{\operatorname{cosk}^{#1}}
\newcommand{\diag}{\operatorname{diag}}
\newcommand{\Der}{\operatorname{Der}}
\newcommand{\Ext}{\operatorname{Ext}}
\newcommand{\Fib}{\operatorname{Fib}}
\newcommand{\gr}{\operatorname{gr}}
\newcommand{\holim}{\operatorname{holim}}
\newcommand{\Hom}{\operatorname{Hom}}
\newcommand{\Id}{\operatorname{Id}}
\newcommand{\Image}{\operatorname{Im}}
\newcommand{\Ker}{\operatorname{Ker}}
\newcommand{\proj}{\operatorname{proj}}
\newcommand{\sk}[1]{\operatorname{sk}_{#1}}
\newcommand{\Sq}[1]{\operatorname{Sq}^{#1}}
\newcommand{\tr}[1]{\operatorname{tr}_{#1}}
\newcommand{\Tot}{\operatorname{Tot}}
\newcommand{\cC}[1]{c^{\lra{#1}}\C}
%
%
\newcommand{\A}{{\mathcal A}}
\newcommand{\Ab}{{\EuScript Ab}}
\newcommand{\Alg}{{\EuScript Alg}}
\newcommand{\C}{{\mathcal C}}
\newcommand{\Ca}{{\mathcal C\hspace{-1mm}A}}
\newcommand{\CA}[1]{\Ca_{#1}}
\newcommand{\CAlg}{{\EuScript{}coAlg}}
\newcommand{\CM}[1]{{#1}\text{-}{\EuScript{}coMod}}
\newcommand{\ca}{co\text{-}ab}
\newcommand{\D}{{\mathcal D}}
\newcommand{\F}{{\mathcal F}}
\newcommand{\K}[1]{{\mathcal K}_{#1}}
\newcommand{\LL}{{\mathcal L}}
\newcommand{\M}{{\mathcal M}}
\newcommand{\cla}{c_{\ast}}
\newcommand{\cua}{c^{\ast}}
\newcommand{\PP}{{\mathcal P}}
\newcommand{\Pu}[1]{{{\EuScript P}}^{#1}}
\newcommand{\Pa}{$\Pi$-algebra}
\newcommand{\RM}{R\text{-}{\EuScript{}Mod}}
\newcommand{\MA}[1]{{\EuScript Mod}\text{-}\A_{p}}
\newcommand{\Set}{{\EuScript Set}}
\newcommand{\Ss}{{\mathcal S}}
\newcommand{\Sa}{\Ss_{\ast}}
\newcommand{\TT}{{\mathcal T}}
\newcommand{\Ta}{\TT_{\ast}}
\newcommand{\U}{{\mathcal U}}
\newcommand{\V}{{\mathcal V}}
%
%
\newcommand{\bD}{\boldsymbol{\Delta}}
\newcommand{\Fp}{{\mathbb F}_{p}}
\newcommand{\Ft}{{\mathbb F}_{2}}
\newcommand{\N}{\mathbb N}
\newcommand{\Q}{\mathbb Q}
\newcommand{\Z}{\mathbb Z}
%
%
\newcommand{\bK}{\mathbf{K}}
\newcommand{\EM}[2]{\bK({#1},{#2})}
\newcommand{\Kp}[2]{{\EuScript E\!M}(R,{#1})_{#2}}
\newcommand{\SKp}[3]{{\EuScript E\!M}^{#1}(R,{#2})_{#3}}
\newcommand{\bS}[1]{\mathbf{S}^{#1}}
\newcommand{\bd}{\mathbf{d}^{0}}
\newcommand{\bbd}{\bar{\mathbf{d}}^{0}}
\newcommand{\dd}[1]{\bd_{#1}}
\newcommand{\bld}{\mathbf{d}_{0}^{\#}}
\newcommand{\bdi}{(\bd)_{\#}}
\newcommand{\ddi}[1]{(\dd{#1})_{\#}}
\newcommand{\qu}[2]{(q^{#1}_{#2})_{\#}}
\newcommand{\bP}[1]{P^{#1}}
\newcommand{\bQ}{\mathbf{Q}}
\newcommand{\bT}{\mathbf{T}}
\newcommand{\W}{\mathbf{W}}
\newcommand{\X}{\mathbf{X}}
\newcommand{\Y}{\mathbf{Y}}
%
%
\newcommand{\co}[1]{c({#1})^{\bullet}}
\newcommand{\Au}{A^{\bullet}}
\newcommand{\Aus}{\Au_{\ast}}
\newcommand{\bAs}[1]{\bar{A}_{\ast}^{#1}}
\newcommand{\Cd}{C_{\bullet}}
\newcommand{\Cud}{C^{\bullet}}
\newcommand{\bC}[1]{\bar{C}^{#1}}
\newcommand{\Du}{\bD^{\bullet}}
\newcommand{\Gus}{G^{\bullet}_{\ast}}
\newcommand{\Qu}{\bQ^{\bullet}}
\newcommand{\baQ}[1]{\bar{\bQ}^{#1}}
\newcommand{\Tu}{\bT^{\bullet}}
\newcommand{\Td}{T_{\bullet}}
\newcommand{\Xd}{X_{\bullet}}
\newcommand{\Xds}{X_{\bullet}^{\ast}}
\newcommand{\Wu}{\W^{\bullet}}
\newcommand{\Xu}{\X^{\bullet}}
\newcommand{\pXu}{X^{\bullet}}
\newcommand{\Yu}{\Y^{\bullet}}
\newcommand{\Yds}{Y_{\bullet}^{\ast}}
%
%
\newcommand{\As}{A_{\ast}}
\newcommand{\Cs}{C_{\ast}}
\newcommand{\Cus}{C^{\ast}}
\newcommand{\Css}{C^{\ast}_{\ast}}
\newcommand{\Gs}{G_{\ast}}
\newcommand{\Js}{J_{\ast}}
\newcommand{\Ks}{K_{\ast}}
\newcommand{\Ls}{L_{\ast}}
\newcommand{\Ms}{M_{\ast}}
\newcommand{\Ts}{T_{\ast}}
\newcommand{\Vs}{V_{\ast}}
\newcommand{\Ws}{W_{\ast}}
\newcommand{\Xs}{X_{\ast}}
\newcommand{\Ys}{Y_{\ast}}
%
%
\newcommand{\Hi}[3]{H_{#1}({#2};{#3})}
\newcommand{\tHi}[3]{\tilde{H}_{#1}({#2};{#3})}
\newcommand{\Hu}[3]{H^{#1}({#2};{#3})}
\newcommand{\His}[2]{\Hi{\ast}{#1}{#2}}
\newcommand{\HR}[2]{\Hi{#1}{#2}{R}}
\newcommand{\HRs}[1]{\His{#1}{R}}
\newcommand{\pus}{\pi^{\ast}}
%
%
\title{Realizing coalgebras over the Steenrod algebra}
\author{David Blanc}
\address{Dept.\ of Mathematics, University of Haifa, 31905 Haifa, Israel}
\email{blanc@math.haifa.ac.il}
\date{Revised: \ March 1, 1999}
\begin{abstract}
We describe algebraic obstruction theories for realizing an abstract 
(co)al\-geb\-ra \ $\Ks$ \ over the mod $p$ Steenrod algebra as the (co)\-homology 
of a topological space, and for distinguishing between the $p$-homotopy types of 
different realizations. The theories are expressed in terms of the Quillen 
cohomology of \ $\Ks$.
\end{abstract}

\maketitle
%
%
\sect{Introduction}
\label{ci}

The question of which graded $R$-algebras can occur as the cohomology
ring of a space $\X$ with coefficients in $R$ was first raised explicitly by
Steenrod in \cite{SteCA}, but it goes back to Hopf, for \ $R=\Q$ \ -- \ see 
\cite{HopfT}. \ When \ $R=\Fp$, \ the cohomology \ $\Hu{\ast}{\X}{\Fp}$ \ 
also has a compatible action of the Steenrod algebra, so it is natural to \vsm 
ask:

Which algebras over the Steenrod algebra can be realized as the cohomology of
a space, and in how many different ways\vsm ?

This question has been addressed repeatedly in the past \ -- \ see, for 
example, \cite{AdHI,AWiF,AdemI,ABNotH1,CEwiR,DMWilH,DWilF,LSmiRS,SSwitR} and 
the survey in \cite{AguRC}. Two related algebraic questions have also often been 
considered: which \ $\Fp$-algebras can be provided with a compatible action of 
the Steenrod algebra (see, e.g., \cite{DKWinsC,STodS,EThomS2}), and conversely, 
which unstable modules over the Steenrod algebra can be provided with a 
compatible algebra structure, or directly: which unstable modules are 
realizable \ -- \ see, e.g., \cite{CSmiRM,KuhnTR,SchwN}.  
The analogous \emph{stable} question of whether a given module over the Steentrod 
algebra can be realized by a spectrum, which has also been extensively studied 
(e.g., \cite{BMadR,BGitS,BPetS}), can be answered in terms of suitable \ $\Ext$ \ 
groups (see \cite[Ch.\ 16, 3]{MargS}), and an unstable version of this for the
Massey-Peterson case was developed by Harper (see \cite{HarpC1}).

However, we shall not be concerned with these variants here: our goal is to 
describe a general obstruction theory for the original realization problem, 
which can be stated purely algebraically, in terms of the Quillen cohomology of 
the given algebra \ -- \ analogous to the stable theory. This answers a question 
of Lannes, Miller and others in the 1980's, asking for an unstable analogue of 
the stable obstruction theory, which was also (independently) one of the 
motivations for the project begun by Dwyer, Kan and Stover in \cite{DKStE,DKStB} 
(see also \cite{BGoeC}).

It turns out to be more natural to consider of the dual question, that of 
realizing a \emph{coalgebra} over the Steenrod algebra as the \emph{homology}
of a space. This is because the cohomology of a space in general has the 
structure of a \emph{profinite} unstable algebra; it is only when the space is 
of finite type that it is an unstable algebra, and the realizability of 
an unstable algebra of finite type is of course strictly equivalent to 
that of the corresponding coalgebra (its vector space dual).

Part the theory we describe here actually works over a more general ground ring 
$R$, but to obtain its full force we restrict attention to the case when $R$ is a 
field. Thus, if we define an \emph{unstable $R$-coalgebra} to 
be a graded coalgebra over $R$ equipped with a compatible action of the 
unstable $R$-homology operations, we have:
%
%
\begin{tha}
For \ $R=\Fp$ \ or $\Q$, \ let \ $\Ks$ \ be a connected unstable $R$-coalgebra, 
such that either \ $K_{1}=0$, \ or \ $\Ks$ \  has finite projective dimension. \ 
Then there is a sequence of cohomology classes \ 
$\chi_{n}\in H^{n+2}(\Ks;\Sigma^{n}\Ks)$ \ such that \ 
$\chi_{n}$ \ is defined whenever \ $\chi_{1}=\dotsc=\chi_{n-1}=0$, \ and 
all the classes vanish if and only if \ $\Ks\cong\HRs{\X}$ \ for some
space $\X$.
\end{tha}

\noindent (See Theorems \ref{tone} and \ref{ttwo} below). The proof involves
showing that, for any space $\X$, any (algebraic) cosimplicial resolution of the 
unstable coalgebra \ $\HRs{\X}$ \ can be realized as a cosimplicial space, and
conversely. As a side benefit, when \ $R=\Fp$, \ this provides a way of 
constructing minimal ``unstable Adams resolutions'' (see Remark \ref{ruar} below).

There is a similar theory for distinguishing between different realizations:
%
%
\begin{thb}
For any two simply-connected spaces $\X$ and $\Y$ of finite type such that \ 
$\His{\X}{\Fp}\cong\His{\Y}{\Fp}\cong\Ks$ \ as unstable coalgebras over the 
Steenrod algebra, there is a sequence of cohomology classes \ 
$\delta_{n}\in H^{n+1}(\Ks;\Sigma^{n}\Ks)$ \ whose vanishing implies that $\X$ 
and $\Y$ are $p$-equivalent.
\end{thb}

\noindent (See Theorems \ref{tthree} and \ref{tfour} below). This result is 
less satisfactory, in that $\X$ and $\Y$ may be $p$-equivalent even if the
cohomology classes do not all vanish (which simply expresses the fact that the 
homotopy theory of cosimplicial spaces is richer than that of topological 
spaces). Note that in this the classes \ $\delta_{n}$ \ resemble the usual 
$k$-invariants, and indeed Theorem B can be thought of as 
providing a system of algebraic ``invariants'' for the $p$-type of a space, 
dual to those of \cite{BlaAI} or \cite{BGoeC} (see \S \ref{rdpt}).

Theorem B also holds for \ $R=\Q$; \ in this case we simply recover 
the homology version of the obstruction theory of \cite{HStaO} and \cite{FelDT}.

\subsection{Notation and conventions}
\label{snac}\stepcounter{thm}

$\TT$ \ will denote the category of topological spaces, and \ $\Ta$ \ that of 
pointed connected topological spaces with base-point preserving maps. We denote
objects in \ $\TT$ or \ $\Ta$ \ by boldface letters: \ $\X,\Y$, \ and so on,
to help distinguish them from the various algebraic objects we consider.

$\Q$ denotes the rationals, and for $p$ prime, \ $\Fp$ \ denotes the field with 
$p$ elements. \ 
For a ring $R$ (always assumed to be commutative with unit), \ $\RM$ \ denotes 
the category of $R$-modules, and \ $R\lra{\Xs}$ \ the free $R$-module on a 
(possibly graded) set of generators \ $\Xs$. \ 
Tensor products of $R$-modules will always be over the ground ring $R$,
unless otherwise stated, and the \emph{dual} module of \ $A\in\RM$ \ is denoted 
by \ $A^{\star}:=\Hom_{\RM}(A,R)$.

$\HRs{\X}\in\RM$ \ is the homology of a topological space (or simplicial set) 
$\X$ with coefficients in $R$. \ We write \ $f_{\#}:\HRs{\X}\to\HRs{\Y}$ \ for
the graded homomorphism induced by \ $f:\X\to\Y$. \ $R_{\infty}\X$ \ is the
Bousfield-Kan $R$-completion of $\X$ \ (cf.\ \cite[I, 4.2]{BKaH}).

For an abelian category $\M$, we let \ $\cla\M$ \ denote the category of chain 
complexes over $\M$ (in non-negative degrees); \ similarly, \ 
$\cua\M$ \ denotes the category of cochain complexes.
 
For any category $\C$, \ we denote by \ $\gr\C$ \ the category of 
\emph{non-\-negatively graded objects} over $\C$, \ with \ 
$\lrv{x}=n\EQUIV x\in X_{n}$ \ for \ $\Xs\in\gr\C$. \ Given a (fixed) object \ 
$B\in\C$, \ we denote by \ $\C\backslash B$ \ the category of objects under $B$ 
(cf.\ \cite[II, \S 6]{MacC}).

\begin{defn}\label{dcos}\stepcounter{subsection}
A \emph{cosimplicial object} \ $\pXu$ \ over any category
$\C$ is a sequence of objects \ $X^{0},X^{1},\ldots,X^{n},\ldots$ \ in
$\C$ equipped with \emph{coface} and \emph{codegeneracy} maps \ 
$d^{i}:X^{n}\to X^{n+1}$, \ $s^{j}:X^{n+1}\to X^{n}$ \ ($0\leq i,j\leq n$) \ 
satisfying the cosimplicial identities
\setcounter{equation}{\value{thm}}\stepcounter{subsection}
\begin{equation}\label{ezero}
\begin{split}
d^{j}d^{i} & = d^{i+1}d^{j}\hsp \text{if} \ i\geq j\\
s^{j}d^{i} & =
\begin{cases} 
d^{i}s^{j-1} &\hsm \text{if $i<j$}\\
\Id           &\hsm \text{if $i=j,j+1$}\\
d^{i-1}s^{j} &\hsm \text{if $i\geq j+2$}\\
 \end{cases}\\
s^{j}s^{i} & = s^{i}s^{j+1}\hsp \text{if} \ i\leq j
\end{split}
\end{equation}
\setcounter{thm}{\value{equation}}
\noindent (cf.\ \cite[X, \S 2.1]{BKaH}). 

We denote by \ $c\C$ \ the category of 
cosimplicial objects over $\C$. \ If we restrict attention to \ 
$X^{0},X^{1},\ldots,X^{n}$, \  with their coface and codegeneracy maps, we 
have an $n$-\emph{cosimplicial} object; the category of such will be denoted by \ 
$\cC{n}$.
\end{defn}

Dually, we denote by \ $s\C$ \ the category of {\em simplicial\/} objects over 
$\C$ \ (cf.\ \cite[\S 2]{MayS}). \ 
The category of simplicial sets, however, will be denoted simply by $\Ss$ 
(rather than \ $s\Set$), \ and that of \emph{pointed} simplicial sets by \ 
$\Sa$. \ Objects in these two categories will again be denoted by boldface 
letters. The standard $n$ simplex in $\Ss$ is denoted by \ 
$\Delta[n]$, \ generated by \ $\sigma_{n}\in \Delta[n]_{n}$, \ and \ 
$\Lambda^{k}_{n}\in\Ss$ \ is the sub-simplicial set of \ $\Delta[n]$ \ 
generated by \ $d_{i}\sigma_{n}$ \ for \ $i\neq k$.

Since we shall be dealing for the most part with simplicial sets as our model for
the homotopy category of topological spaces, we shall call cosimplicial 
pointed simplicial sets \ -- \ i.e., objects in \ $c\Sa$ \ -- \ simply 
\emph{cosimplicial spaces}. 

\begin{example}\label{ecs}\stepcounter{subsection}
The cosimplicial space \ $\Du\in c\Ss$ \ has the standard simplicial 
$n$-simplex \ $\Delta[n]$ \ in cosimplicial dimension $n$, with coface and 
codegeneracy maps being the standard inclusions and projections (cf.\ 
\cite[I, 3.2]{BKaH}).
\end{example}

\begin{defn}\label{dcs}\stepcounter{subsection}
If $\C$ has enough limits, the obvious truncation functor \ 
$\tr{n}:c\C\to\cC{n}$ \ has a right adjoint \ $\rho_{n}$, \ and the 
composite \ $\cosk{n}\DEF\rho_{n}\circ\tr{n}:c\C\to c\C$ \ is called the 
$n$-\emph{coskeleton} functor. (This is dual to $n$-\emph{skeleton} functor \ 
$\sk{n}:s\C\to s\C$.)
\end{defn}

\subsection{Organization}
\label{sorg}\stepcounter{thm}

In section \ref{cuc} we recall some basic facts about coalgebras over the 
Steenrod algebra, and in section \ref{crc} we show how certain convenient CW 
resolutions for such coalgebras may be constructed.
Section \ref{cfc} deals with the coaction of the fundamental group of a
cosimplicial coalgebra, and the Quillen cohomology of unstable coalgebras.
In section \ref{crr} we describe the cohomology classes which determine whether 
a given algebraic resolution may be realized topologically (Theorem \ref{tone}),
and in section \ref{crec} we apply this to the original question of realizing an 
abstract coalgebra (Theorem \ref{ttwo}). Finally, in section \ref{cdr} a similar
theory is developed for distinguishing
between different possible realizations of a given algebraic resolution 
(Theorem \ref{tthree}), and thus for determining the $p$-type of a space 
(Theorem \ref{tfour}).

\subsection{Acknowledgements}
\label{sack}\stepcounter{thm}
I would like to thank Pete Bousfield, Bill Dwyer, and Brooke Shipley for several
useful conversations, and Dan Kan for suggesting that I consider the problem of 
realizing algebras over the Steenrod algebra, which was one of the 
motivations for the project begun in \cite{DKStE} and \cite{DKStB}. \ 
I would also like to thank the referee for many helpful comments and suggestions.

%
%
\sect{Unstable coalgebras}
\label{cuc}

We first recall some basic facts about the category of the coalgebras over the 
Steenrod algebra:

\begin{defn}\label{dca}\stepcounter{subsection}
For a field $R$, let \ $\CAlg_{R}$ \ denote the category of graded coalgebras 
over $R$: an object in \ $\CAlg_{R}$, \ which we shall call simply an 
\emph{$R$-coalgebra}, is thus a (non-negatively) graded $R$-module \ 
$\Vs\in\gr\RM$, \ equipped with a coassociative \emph{diagonal} 
(or \emph{comultiplication}) map \ $\Delta:\Vs\to\Vs\otimes\Vs$, \ 
with \ $(\Id\otimes\Delta)\circ \Delta=(\Delta\otimes\Id)\circ \Delta$, \ 
and an \emph{augmentation} (or \emph{counit}) map \ $\varepsilon:\Vs\to R$, \ 
with \ $(\Id\otimes\varepsilon)\circ \Delta$ \ and \ 
$(\varepsilon\otimes\Id)\circ\Delta$ \ equal respectively to the natural 
isomorphisms \ $\Vs\to\Vs\otimes R$ \ and \ $\Vs\to R\otimes\Vs$. \ 
We require the comultiplication to be \emph{cocommutative}, in the graded 
sense \ -- \ i.e., \ $\Delta\circ\tau=\Delta$, \ where \ 
$\tau(a\otimes b)\DEF -1^{\lrv{a}\lrv{b}}b\otimes a$ \ is the graded switch 
map. \ See \cite[\S 1.0]{SweeH} and \ \cite[\S 2.1]{MMoorH}.

We assume all our graded coalgebras \ $\Cs\in\CAlg_{R}$ \ are \ 
\emph{connected} \ -- that is, \ $C_{0}\cong R$; \ $\Cs\in\CAlg_{R}$ \ is called 
\emph{simply-connected} if in addition \ $C_{1}=0$. \ The coalgebra \ $\Cs$ \ is 
of \emph{finite type} if \ $C_{k}$ \ is finite dimensional vector space over $R$ 
for each \ $k\geq 0$. \ We 
can pass from \ $\CAlg_{R}$ \ to the category \ 
$\Alg_{R}$ \ of connected, unital, graded-commutative \emph{algebras} over $R$ 
by taking the vector-space dual: \ $(\Cs)^{\star}\in\Alg_{R}$ \ (and if \ 
$\Cus\in\Alg_{R}$ \ is of finite type, we can of course pass back to \ 
$\CAlg_{R}$ \ in the same way).
\end{defn}

\subsection{Unstable coalgebras}
\label{suc}\stepcounter{thm}
 
As usual, the homology \ $\HRs{\X}$ \ of a space $\X$ with coefficients in field 
$R$ is an $R$-coalgebra, with \ $\Delta$ induced by the diagonal \ 
$\X\to\X\times\X$ \ (composed with inverse of the K\"{u}nneth isomorphism \ 
$\HRs{\X}\otimes\HRs{\X}\cong\HRs{\X\times\X}$). \ 
However, \ $\HRs{\X}$ \ also comes equipped with an action of the primary 
$R$-homology operations: \ these are natural transformations \ 
$\HR{i}{-}\to\HR{k}{-}$, \  dual to the corresponding cohomology operations, and 
they vanish if \ $k>i$ \ (see \cite[\S 9]{SteCOP}).

For any field $R$, an \emph{unstable coalgebra} (over $R$) is a non-negatively
graded $R$-module equipped with an action of primary $R$-homology operations
(which include the coalgebra structure), satisfying the universal identities
for these operations. We denote the category of such unstable $R$-coalgebras by \ 
$\CA{R}$.

The simplest case is when \ $R=\Q$: \ $\CA{\Q}\approx \CAlg_{\Q}$, \ since there 
are no non-trivial primary $\Q$-homology operations besides the coproduct 
(see \cite[I, \S 1]{QuR}). The next simplest is \ $R=\Fp$:

\begin{defn}\label{dum}\stepcounter{subsection}
For any prime $p$, an \emph{unstable module over the mod $p$ Steenrod algebra}, \ 
$\A_{p}$, \ is a non-negatively graded \ $\Fp$-vector space \ $\Ks$, \ equipped 
with a right action of $\A_{p}$ \ -- \ i.e., a graded homomorphism \ 
\setcounter{equation}{\value{thm}}\stepcounter{subsection}
\begin{equation}\label{eone}
\lambda:\Ks\otimes\A_{p}\to\Ks
\end{equation}
\setcounter{thm}{\value{equation}}
\noindent where $\lrv{x\Sq{i}}=\lrv{x}-i$ \ if \ $p=2$, \ and \ 
$\lrv{x\Pu{i}}=\lrv{x}-2(p-1)i$, \ $\lrv{x\beta}=\lrv{x}-1$ \ if \ 
$p>2$ \ (where we write \ $x\Sq{i}$ \ for \ $\lambda(x\otimes\Sq{i})$, \ etc.). 
The action is required to be \emph{unstable} in the sense that \ 
$x\Sq{i}=0$ \ if \ $2i>\lrv{x}$ \ (for \ $p=2$) \ and \ 
$x\Pu{i}=0$ \ if \ $2pi>\lrv{x}$ \ (for \ $p>2$). \ 
The category of such unstable modules will be denoted by \ $\MA{\Fp}$.
\end{defn}

The category \ $\MA{\Fp}$ \ is dual to the more familiar category $\U$ of 
unstable ``coho\-mo\-lo\-gy-\-like'' modules over the Steenrod algebra (see, 
e.g., \cite[\S 1.3]{SchwU}).

\begin{defn}\label{dua}\stepcounter{subsection}
For any prime $p$, the category \ $\CA{\Fp}$ \ of 
graded coalgebras over the mod $p$ Steenrod algebra \ $\A_{p}$, \ has as objects
non-negatively graded coalgebras \ $\Cs$ \ over \ $\Fp$, \ which are at the 
same time unstable \ $\A_{p}$-modules. The two structures are related by the 
Cartan formula, which says that the $\lambda$ of \eqref{eone} is a homomorphism
of coalgebras (see \cite[\S 4]{MilnS}) \ -- \ dual to: \ 
$\Sq{n}(a\cdot b)=\sum_{k}\Sq{k}a\cdot\Sq{n-k}b$ \ for \ $p=2$, \ and \ 
$\Pu{n}(a\cdot b)=\sum_{k}\Pu{k}a\cdot\Pu{n-k}b$, \ $\beta(a\cdot b)=
(\beta a)\cdot b+(-1)^{\lrv{a}}a\cdot\beta b$ \ for \ $p>2$. \ 
There is also a Verschiebung formula, dual to fact that that the top Steenrod
operation equals the Frobenius \ -- \ i.e., for \ $\lrv{a}=n$, \ we have \ 
$\Sq{n}a=a^{2}$ \ if \ $p=2$, \ and \ $\Pu{n/2}a=a^{p}$ \ if \ $p>2$ \ and $n$ 
is even. \ See \cite[\S 5]{BCurS}.
\end{defn}

In particular, we can think of\ $\His{\X}{\Fp}$ \ as an object in either \ 
$\MA{\Fp}$ \ or \ $\CA{\Fp}$ \ for any space $\X$.
Again, \ $\CA{\Fp}$ \ is dual to the more familiar category \ $\K{}=\K{\Fp}$ \ of 
unstable algebras over the Steenrod algebra (cf.\ \cite[\S 1.4]{SchwU}). \ 
However, taking vector space duals yields a strict equivalence of categories
only when dealing with (co)algebras of \emph{finite type}, and our approach is 
more naturally presented in terms of coalgebras, as noted above.

\subsection{(Co)abelian (co)algebras}
\label{sac}\stepcounter{thm}

Recall that  any \emph{abelian} $R$-algebra (i.e., abelian group object in \ 
$\Alg_{R}$) \ must have a trivial multiplication (for any ring $R$). \ 
When  \ $R=\Q$, \ the subcategory \ $(\K{\Q})_{ab}$ \ of 
abelian objects in \ $\K{\Q}$ \ is actually equivalent to the 
category of graded vector spaces over $\Q$. \ $(\K{\Fp})_{ab}$ \ is equivalent to 
a subcategory of \ $\MA{\Fp}$ \ (viewed as algebras with a trivial product):\hs \ 
for \ $p=2$, \ $(\K{\Ft})_{ab}=\Sigma\U$ \ is the category of 
$\A_{2}$-modules with \ $Sq^{i}x=0$ \ for \ $|x|\leq i$; \ \ for \ $p>2$, \ 
$(\K{\Fp})_{ab}=\V$ \ is the category of $\A_{p}$-modules with \ $\PP^{i}x=0$ \ 
for \ $|x|\leq 2i$ \ (cf.\ \cite[\S 1]{MilSC}). \ In all these cases the
abelianization functor \ $()_{ab}:\K{R}\to(\K{R})_{ab}$ \ assigns to any 
algebra \ $A^{\ast}\in\K{R}$ \ its ``module of indecomposables'', \ $Q(A^{\ast})$.

Dualizing, we see that the \emph{coabelian} objects (i.e., abelian cogroup 
objects) in \ $\CA{R}$ \ must have trivial comultiplication, so \ 
$(\CA{R})_{\ca}$ \ is equivalent to the appropriate subcategory of \ $\MA{R}$; \ 
the coabelianization functor is just the $R$-module of primitives \ 
$\PP(\As)\hra \As$, \ for any \ $\As\in\CA{R}$ \ (see \cite[\S 8.6]{BousH}, 
and compare \cite[\S 2]{KLeeC}).

\subsection{Functors and limits of coalgebras}
\label{sflc}\stepcounter{thm} 

The underlying-set functor \ $\CAlg_{R}\to\gr\Set$ \ factors through \ 
$\hat{U}:\CAlg_{R}\to\gr\RM$, \ with right adjoint \ 
$\hat{G}:gr\RM\to\CAlg_{R}$, \ where \ $\hat{G}(\Vs)$ \ is the (cocommutative) 
\emph{cofree coalgebra} on \ $\Vs$ \ (cf.\ \cite[\S 6.4]{SweeH}). 
Moreover, the functor $\hat{G}$ creates all colimits 
in \ $\CAlg_{R}$, \ in the sense of \cite[V, \S 1]{MacC}, and the pair \ 
$(\hat{U},\hat{G})$ \ \emph{produces} all limits in \ $\CAlg_{R}$ \ in the sense
of \cite[\S 3.3]{BlaN}. \ The same is true for the right adjoint \ 
$G:gr\RM\to\CA{R}$ \ of the ``underlying graded $R$-module'' functor \ 
$U:\CA{R}\to\gr\RM$. \ See \cite[Prop.\ 7.5]{BlaN} and \cite[\S 8.2]{BousH}. 

Note also that any $R$-coalgebra, as well as any unstable $R$-coalgebra, is 
isomorphic to the colimit of its finite sub-coalgebras (see \cite[1.1]{GoeHH}),
and this allows one to describe the product of an arbitrary collection of
coalgebras \ $(\Cs^{(i)})_{i\in I}$ \ in terms of the partially
ordered collection $\F$ of finite subsets \ $J\subseteq I$ \ as \ 
\setcounter{equation}{\value{thm}}\stepcounter{subsection}
\begin{equation}\label{etwo}
\prod_{i\in I}\Cs^{(i)}=\colim_{J\in\F}\colim_{\alpha \in A_{J}}
(\otimes_{i\in J}\Cs^{(i)})_{\alpha},
\end{equation}
\setcounter{thm}{\value{equation}}
\noindent where \ $(\otimes_{i\in J}\Cs^{(i)})_{\alpha}$ \ ($\alpha\in A_{J}$) \ 
runs over all \emph{finite} sub-coalgebras of the finite tensor product \ 
$\otimes_{i\in J}\Cs^{(i)}$ \ (see \cite[1.2]{GoeHH}).

\begin{remark}\label{rcuc}\stepcounter{subsection}
In fact, any \emph{cofree unstable $R$-coalgebra}  \ $\Gs=G(\Vs)\in\CA{R}$ \ is 
of the form \ $\Gs\cong\HRs{\bK(R\lra{\Xs})}$, \ where \ $\Xs\in\gr\Set$ \ 
is a graded set and \ $\bK(R\lra{\Xs})$ \ is the GEM (generalized 
Eilenberg-Mac Lane object) \ $\prod_{n=1}^{\infty}\prod_{x\in X_{n}}\EM{R}{n}$. \ 
By \eqref{etwo} we have \ $\Gs\cong \prod_{n=1}\prod_{x\in X_{n}}\HRs{\EM{R}{n}}$.

Note that each \ $\HRs{\EM{R}{n}}$, \ and thus \ $\Gs$, \ is an abelian Hopf 
algebra (see \cite[\S 4.4]{BousHS}), \ and for any map of graded sets \ 
$f:\Xs\to\Ys$, \ the map \ $G(f):G(\Xs)\to G(\Ys)$ \ is a morphism of Hopf 
algebras (in particular, an algebra homomorphism).
\end{remark}
%
%
\sect{Resolutions of coalgebras}
\label{crc}

We now prove some basic facts about cofree resolutions for coalgebras:

\begin{defn}\label{dmat}\stepcounter{subsection}
For any concrete cocomplete category $\C$, the \emph{co-matching object} 
functor \ $M:\Ss^{op}\times c\C\to \C$, \ written \ $M^{K}\pXu$ \ for \ a finite 
simplicial set \ $K\in\Ss$ \ and  any \ $\pXu\in c\C$, \ is defined by
requiring that \ $M^{\Delta[n]}\pXu:=X^{n}$, \ and if \ $K=\colim_{i} K_{i}$ \ 
then \ $M^{K}\pXu=\colim_{i} M^{K_{i}}\pXu$ \ (cf.\ \cite[\S 6]{BousHS} and 
\cite[\S 2.1]{DKStB}).

In particular, write \ $M^{n}\pXu$ \ for \ $M^{\sk{n-1}\Delta[n]}\pXu$. \ 
Note that each coface map \ $d^{k}:X^{n-1}\to X^{n}$ \ factors through the map \ 
$\xi^{n}:M^{n}\pXu\to X^{n}$ \ induced by the inclusion \ 
$\sk{n-1}\Delta[n]\hra\Delta[n]$. \ 
A cosimplicial space \ $\Xu\in c\Sa$ \ will be called \emph{cofibrant} if each
of these maps \ $\xi^{n}$ \ is a cofibration. This concept refers to the 
resolution model category structure on \ $c\Sa$ \ (see \S \ref{srmc} below), 
rather than the Reedy model category structure of \cite[X,\S 4]{BKaH}.
\end{defn}

\begin{defn}\label{dlat}\stepcounter{subsection}
Given a cosimplicial object \ $\pXu$ \ over a concrete complete category $\C$, \ 
the analogous construction for the codegeneracies yields \ $L^{n}\pXu$ \ 
is defined (in the cases of interest to us) by \ 
$$
L^{n}\pXu\DEF\{(x_{0},\dotsc,x_{n-1})\in(X^{n-1})^{\times n}~\lvert\ 
s^{i}x_{j}=s^{j-1}x_{i}\ \ \text{for all}\ 0\leq i<j\leq n\}.
$$
\noindent Again, each codegeneracy map \ $s^{i}:X^{n}\to X^{n-1}$ \ equals the 
natural map \ $\vars^{n}:X^{n}\to L^{n}\pXu$, \ composed with the projection onto
the $i$-th factor. 

$L^{n}\pXu$  has been called the $n$-th ``co-latching object'' for \ 
$\pXu$ \ -- \ cf.\ \cite[\S 2.3]{DKStE}. In \cite[X, \S 4.5]{BKaH} it is
denoted by \ $M^{n}\Xu$; \ the notation we have here was chosen to be consistent 
with that of \cite{DKStE,DKStB} and \cite{BlaAI}.
\end{defn}

\begin{defn}\label{dnc}\stepcounter{subsection}
If \ $\Xu\in c\Sa$ \ is cofibrant, \ its $n$-\emph{cochains} object,  written \ 
$C^{n}\Xu$, \ is defined to be the cofiber of \ 
$\xi^{n}_{0}:M^{\Lambda^{0}_{n}}\Xu\to\X^{n}$ \ (\S \ref{dcos}), \ 
so \ $C^{n}\Xu=\X^{n}/(\cup_{i=1}^{n}d^{i}\X^{n-1})$. \ 
Similarly, the object \ $B^{n}\Xu$ \ is defined to be the cofiber of \ 
$\xi^{n}:M^{n}\Xu\to\X^{n}$, \ so that \ 
$B^{n}\Xu=\X^{n}/(\cup_{i=0}^{n}d^{i}\X^{n-1})$.
\end{defn}

These all fit into the commutative diagram of Figure \ref{fig1}:
\newpage
%
%
\begin{figure}[htbp]
\begin{picture}(180,65)(-120,0)
%
%
\put(25,50){$\X^{n}$}
\put(45,55){\vector(1,0){50}}
\put(65,58){$d^{0}$}
\put(100,50){$\X^{n+1}$}
%
%
\put(35,48){\vector(0,-1){33}}
\put(35,15){\vector(0,-1){2}}
\put(-5,30){$q^{n}_{X}\circ p^{n}_{X}$}
\put(115,48){\vector(0,-1){33}}
\put(115,15){\vector(0,-1){2}}
\put(92,30){$p^{n+1}_{X}$}
\put(130,55){\vector(3,-2){60}}
\put(190,15){\vector(3,-2){2}}
\put(162,40){$q^{n+1}_{X}\circ p^{n+1}_{X}$}
%
%
\put(20,0){$B^{n}\Xu$}
\put(55,5){\vector(1,0){2}}
\put(57,5){\vector(1,0){40}}
\put(67,11){$\bd_{\X^{n}}$}
\put(100,0){$C^{n+1}\Xu$}
\put(145,5){\vector(1,0){40}}
\put(185,5){\vector(1,0){2}}
\put(150,12){$q^{n+1}_{X}$}
\put(193,0){$B^{n+1}\Xu$}
\end{picture}
\caption[fig1]{}
\label{fig1}
\end{figure}

\noindent in which \ $p^{n}_{X}:\X^{n}\epic C^{n}\Xu$ \ and \ 
$q^{n}_{X}:C^{n}\Xu\epic B^{n}\Xu$ \ are quotient maps, and the bottom row is a 
cofibration sequence (with the cofibration \ 
$\bd=\bd_{X^{n}}:B^{n}\Xu\to C^{n+1}\Xu$ \ 
induced from \ $\xi^{n}:M^{n}\Xu\to X^{n}$ \ by cobase change).
We will call \ $\bd$ \ the $n$-th \emph{principal face map} for \ $\Xu$.

\begin{defn}\label{dcoh}\stepcounter{subsection}
The same definitions may be applied to a cosimplicial object \ $\Au\in c\C$ \ 
over any suitable category $\C$ \ (e.g., if $\C$ has a faithful 
``underlying object'' functor into an abelian category).
In this case the kernel of \ $\bd_{A}=B^{n}\Au\to C^{n+1}\Au$ \ is defined as 
usual to be the $n$-\emph{th cohomotopy object} of \ $\Au$, \ and denoted by \ 
$\pi^{n}\Au$ \ (see \cite[X, \S 7.1]{BKaH}). \ When $\C$ is an abelian
category, \ $\Au\in c\C$ \ is equivalent under the Dold-Kan equivalence 
(cf.\ \cite[Thm 1.9]{DoH}) to a cochain complex \ $A^{\ast}$, \ and 
$\pi^{n}\Au\cong H^{n}A^{\ast}$. \ 
However, in most cases \ $\pi^{n}\Au$ \ will have additional structure \ -- \ 
e.g., it will be a (coabelian) object in \ $\C$ \ (see \ \cite[\S 5]{BStG}).
\end{defn}

\begin{remark}\label{rcoh}\stepcounter{subsection}
Let \ $\Sigma^{n}\EM{R}{k}\in c\Sa$ \ denote the cosimplicial space consisting of
the usual simplicial Eilenberg-Mac Lane space \ $\EM{R}{k}$ \ in cosimplicial 
dimension $n$, a single point in lower dimensions, and \ 
$L^{p}\Sigma^{n}\EM{R}{k}$ \ (defined inductively) in dimension \ $p>n$; \ 
similarly \ $C\Sigma^{n}\EM{R}{k}$ \ is obtained from \ $\Sigma^{n}\EM{R}{k}$ \ 
by attaching a single copy of \ $\EM{R}{k}$ \ in cosimplicial dimension \ 
$n+1$ \ (see \S \ref{dcw} and compare \cite[\S 3.6]{DKStB}). \ Write \ 
$\SKp{n}{k}{\Ca}$ \ for the cosimplicial unstable coalgebra \ 
$G(R\lra{\Sigma^{n}\EM{R}{k}})$ \ (and similarly \ 
$C\SKp{n}{k}{\Ca}\DEF G(R\lra{C\Sigma^{n}\EM{R}{k}})$). \ Then for any \ 
$\Aus\in c\CA{R}$ \ we have \ 
$\Hom_{c\CA{R}}(\Aus,\SKp{n}{k}{\Ca})\cong B^{n}\Au_{k}$ \ and 
$\Hom_{c\CA{R}}(\Aus,C\SKp{n}{k}{\Ca})\cong C^{n+1}\Au_{k}$, \ so it 
makes sense to denote \ $\pi^{n}\Au_{k}$ \ by \ $[\Aus,\SKp{n}{k}{\Ca}]$ \ -- \ 
and these are in fact the homotopy classes of maps into \ 
$\SKp{n}{k}{\Ca}$ \ in the model category structure on \ $c\CA{R}$ \ 
described in \cite[\S 7]{BlaN}.
\end{remark}

The following statement is dual to \cite[Lemma 2.29]{BlaAI}:
%
%
\begin{prop}\label{pone}\stepcounter{subsection}
For any field $R$ and cofibrant \ $\Xu\in c\Sa$, \ the inclusion \ 
$\iota:C^{n}\Xu\hra\X^{n}$ \ induces an isomorphism \ 
$\iota_{\star}:\HRs{C^{n}\Xu}\cong C^{n}(\HRs{\Xu})$ \ for each \ $n\geq 0$.
\end{prop}

\begin{proof}
The free $R$-module functor \ $F:\Sa\to s\RM$ \ is a left adjoint, so preserves 
all colimits, and thus \ $C_{n}(F\Xu)\cong F(C^{n}\Xu)$ \ (we need \ 
$\Xu$ \ to be cofibrant in order for \ $C^{n}\Xu$ \ to be meaningful). \ 
If \ $K:s\RM\xra{\cong}\cla\RM$ \ is the Dold-Kan equivalence functor, then for 
any \ $\X\in\Sa$, \ $\HRs{\X}$ \ is the homology of the chain complex \ $KF\X$, \ 
so it suffices to show that for any cosimplicial chain complex \ $\Aus(=KF\Xu)$ \ 
the map \ $\iota_{\star}:H_{k}(C^{n}\Aus)\cong C^{n}(H_{k}\Aus)$ \ is an 
isomorphism for all \ $k,n$.

Now for a cosimplicial object \ $\Au$ \ over any abelian category, we can 
define the \emph{Moore cochain complex} by \ 
$N^{n}\Au\DEF \cap_{j=0}^{n-1}\Ker(s^{j})\subseteq A^{n}$ \ (with differential \ 
$\delta\DEF\sum_{i=0}^{n}(-1)^{i}d^{i}$). \ We claim that the composite \ 
$N^{n}\Au\hra A^{n}\epic C^{n}\Au$ \ is an isomorphism \ 
$\Phi:N^{n}\Au\xra{\cong}C^{n}\Au$.

First note that \ $N^{n}\Au\cap\bigcup_{i=1}^{n}\Image(d^{i})=0$, \ since if \ 
$\alpha=\sum_{i=\ell}^{n}d^{i}(x_{i})$ \ (which we may assume by induction on \ 
$0\leq \ell\leq n$) \ and \ $s^{j}\alpha=0$ \ for \ $0\leq j\leq n$, \ then \ 
$0=s^{\ell-1}\alpha=x_{\ell}+\sum_{i=\ell+1}^{n}d^{i-1}s^{\ell-1}(x_{i})$, \ so \ 
$d^{\ell}x_{\ell}=-\sum_{i=\ell+1}^{n}d^{i}d^{\ell}s^{\ell-1}(x_{i})$. \ Thus 
$\alpha$ is in fact of the form \ $\sum_{i=\ell+1}^{n}d^{i}(x'_{i})$ \ -- \ 
which implies that \ $\alpha=0$ \ (for \ $\ell=n$). \ This shows that \ 
$\Phi$ \ is one-to-one.

Next, given \ $\alpha\in A^{n}$ \ with \ $s^{j}(\alpha)=0$ \ for \ 
$0\leq j<\ell$ \ 
(which we again assume by induction on \ $0\leq \ell<n$), \ then there is an \ 
$\alpha'\in A^{n}$ \ with \ $[\alpha]=[\alpha']\in C^{n}\Au$ \ such that \ 
$s^{j}(\alpha')=0$ \ for \ $0\leq j\leq \ell$ \ -- \ namely, \ 
$\alpha'\DEF\alpha-d^{\ell+1}s^{\ell}\alpha$. \ This shows that $\Phi$ is onto.

Thus the Proposition will follow if we show \ 
$\iota_{\star}:H_{k}(N^{n}\Aus)\cong N^{n}(H_{k}\Aus)$ \ is an isomorphism:
given \ $\lra{\alpha}\in N^{n}(H_{k}\Aus)$ \ represented by \ 
$\alpha\in A^{n}_{k}$ \ with \ $\partial_{k}(\alpha)=0$, \ (where $\partial$ 
is the boundary map of the chain complex \ $A^{n}_{\ast}$), \ we assume that \ 
$s^{j}(\alpha)=0$ \ for \ $0\leq j<\ell$, \ and that there are \ 
$b_{\ell},\dotsc,b_{n-1}\in A^{n-1}_{k+1}$ \ such that \ 
$s^{j}(\alpha)=\partial_{k+1}(b_{j})$ \ for \ $\ell\leq j<n$. \ Replacing 
$\alpha$ by \ $\alpha'\DEF \alpha-\partial_{k+1}(d^{\ell+1}(b_{\ell}))$, \ we 
see by induction on  \ $0\leq \ell\leq n-1$ \ that we may choose a representative 
for \ $\lra{\alpha}$ \ with \ $s^{j}(\alpha)=0$ \ for all $j$. \ Thus \ 
$\iota_{\star}$ \ is surjective. Finally, if \ 
$\lra{\alpha}=0\in N^{n}(H_{k}\Aus)$ \ for \ $\alpha\in N^{n}\Au_{k}$, \ there is 
a \ $\beta\in A^{n}_{k+1}$ \ such that \ $\partial_{k+1}(\beta)=\alpha$, \ with \ 
$s^{j}(\beta)=0$ \ for \ $0\leq j\leq\ell$; \ setting \ 
$\beta'\DEF \beta-d^{\ell+1}s^{\ell+1}(\beta)$, \ we see that we can assume \ 
$\beta\in N^{n}B_{k+1}$, \ so \ $\iota_{\star}$ \ is one-to-one.
\end{proof}

\begin{defn}\label{dfsr}\stepcounter{subsection}
A cosimplicial coalgebra \ $\Aus$ \ is called \textit{cofree} if for each \
$n\geq 0$ \ there is a graded set \ $\Ts^{n}$ \ of elements in \ $\As^{n}$ \
such that \ $\As^{n}=G(\Ts^{n})$ \ (cf.\ \S \ref{sflc}), \ and 
\setcounter{equation}{\value{thm}}\stepcounter{subsection}
\begin{equation}\label{ethree}
\text{each codegeneracy map} \ s^{j}:\As^{n}\to \As^{n-1} \ \text{takes} \ 
\Ts^{n} \ \text{to} \ \Ts^{n-1}.
\end{equation}
\setcounter{thm}{\value{equation}}
\end{defn}

\begin{defn}\label{dcw}\stepcounter{subsection}
A \emph{CW complex} over a pointed category $\C$ is a cosimplicial object \ 
$\Cud\in s\C$, \ together with a sequence of objects \ $\bC{n}\in\C$ \ 
($n=0,1,\dotsc$) \ -- called a \emph{CW basis} for \ $\Cud$ \ -- \ such that \ 
$C^{n}=\bC{n}\times L^{n}\Cud$ \ for each \ $n\geq 0$, \ with projection \ 
$\proj_{\bar{C}^{n}}:C^{n}\to\bC{n}$, \ such that \ 
$\proj_{\bar{C}^{n}}\circ d^{i}=0$ \ for \ $1\leq i\leq n$ \ 
(compare \cite[\S 3]{KanR} and \cite[\S 5]{BlaD}).

The coface map \ 
$\bar{d}_{n}^{0}\DEF \proj_{\bar{C}^{n}}\circ d^{0}:C^{n-1}\to\bC{n}$ \ 
is called the \emph{attaching map} for \ $\bC{n}$, \ and it is readily verified 
that the attaching maps \ $\bar{d}_{n}^{0}$ \ ($n=0,1,\dotsc$), \ together with 
the cosimplicial identities \eqref{ezero}, determine all the face and degeneracy 
maps of \ $\Cud$. \ Note that \ $\bar{d}_{n}^{0}$ \ factors through a map \ 
$B^{n-1}\Cud\to\bC{n}$.

In particular, we require that a CW basis for a free cosimplicial algebra \ 
$\Aus\in\CA{R}$ \ be a sequence of \emph{cofree} colagebras \ 
$(\bAs{n})_{n=0}^{\infty}$. 

On the other hand, for a cosimplicial object over \ $\Sa$, \ it will be
convenient to require only that \ 
$C^{n}\simeq\bC{n}\times L^{n}\Cud$ \ for all $n$.
\end{defn}

For any field $R$, the category \ $c\CA{R}$ \ of cosimplicial unstable 
$R$-coalgebras has a model category structure in the sense of Quillen (cf.\ 
\cite[I, \S 1]{QuH}), induced from that of \ $c\RM\approx \cua\RM$ \ 
by the obvious pair of adjoint functors (see \cite[\S 7]{BlaN}). \ All we 
shall need from this are the following:

\begin{defn}\label{dccr}\stepcounter{subsection}
A \textit{cofree cosimplicial resolution} of an unstable $R$-coalgebra  \ 
$\Ks$ \ is defined to be a cofree cosimplicial coalgebra \ 
$\Aus$, \ equipped with a coaugementation \ $\Ks\to\As^{0}$, \ such that 
in each degree \ $k\geq 1$, \ the cohomotopy groups of 
the cosimplicial $R$-module \ $(\Au)_{k}$ \ (i.e., the cohomology groups of the 
corresponding cochain complex over \ $\RM$) \ vanish in dimensions \ 
$n\geq 1$, \ and the coaugmentation induces an isomorphism \ 
$\Ks\cong\pi^{0}\Aus$.
\end{defn}

\subsection{Constructing CW resolutions}
\label{sccr}\stepcounter{thm}

As usual, such a resolution is simply a fibrant (and cofibrant) object in \ 
$c\CA{R}$ \ which is weakly equivalent to the constant cosimplicial object \ 
$\co{\Ks}$, \ in the model category structure mentioned above. In particular,
such resolutions always exist, for any \ $\Ks\in\CA{R}$; \ there are a number
of ways to construct them, including the (very large) canonical monad resolution
described in \cite[\S 11.4]{BKaS} when $R$ is a field (see \cite[\S 7.8]{BlaN}).

We shall be interested in a particular type, namely, those equipped with a 
CW basis, which will be called \emph{CW resolutions}, since these can be chosen
to be small (e.g., minimal). Their construction is straightforward: starting 
with \ $\As^{-1}\DEF\Ks$ \ and \ $B^{-1}\Aus\DEF\As^{-1}$, \ we assume that we
have constructed \ $\Aus$ \ through cosimplicial dimension \ $n-1$; \ then we
simply choose some cofree unstable coalgebra \ $\bAs{n}\in\CA{R}$ \ with 
a one-to-one map \ $\bar{d}^{0}_{n}:B^{n-1}\Aus\hra\bAs{n}$. \  
These always exist, by the universal property of cofree coalgebras \ -- \ e.g., 
one could take \ $\bAs{n}\DEF GU(B^{n-1}\Aus)$. \ Setting \ 
$\As^{n}\DEF\bAs{n}\times L^{n}\Aus$ \ completes the inductive stage. 

The dual construction, for simplicial groups, algebras, and so on, is classical:
see \cite{KanR}, \cite{TateH} and \cite[I, \S 6]{AndrM}. \ 
Note, however, the following analogue of \cite[Prop. 3.18]{BlaC}:

%
%
\begin{prop}\label{ptwo}\stepcounter{subsection}
Any cofree cosimplicial unstable coalgebra \ $\Aus\in c\CA{R}$ \ has a CW 
basis \ $(\bAs{n})_{n=0}^{\infty}$.
\end{prop}

\begin{proof}
Start with \ $\bAs{0}\DEF\As^{0}$, \ and note that \eqref{ethree} of definition 
\ref{dfsr} implies (by induction on $n$) that \ $L^{n}\Aus\cong G(R\lra{\Ys})$ \ 
for some \ $\Ys\in\gr\Set$. \ 

Now because \ $\CA{R}$ \ has the ``underlying structure'' of an abelian category,
we may define a homomorphism of the underlying abelian groups \ 
$\psi^{n}:\As^{n}\to\As^{n}$ \ by \ 
$$
\psi^{n}(\alpha)\DEF \sum_{k=1}^{[(n+1)/2]}\ \sum_{(I,J)\in\LL_{k}^{n}}\ 
(-1)^{\lrv{I}+\lrv{J}+1} d^{I}s^{J}\alpha,
$$
\noindent where \ $\LL_{k}^{n}=\{(I,J)\in\N^{k}\times \N^{k}~\lvert\ 
j_{k-t} >i_{t}> j_{k+1-t}\ \text{for all}\ 1\leq t\leq k\}$ \ (and we let \ 
$j_{0}\DEF n$).

It then follows from the cosimplicial identities \ 
\eqref{ezero}\ that \ $s^{j}\psi^{n}(\alpha)=s^{j}\alpha$ \ for \ 
$0\leq j\leq n-1$. \ Since the definition of \ $\psi^{n}$ \ depends only on \ 
$(s^{0}\alpha,\dotsc,s^{n-1}\alpha)\in L^{n}\Aus$, \ we see that \ 
$\vars^{n}:\As^{n}\to L^{n}\Aus$ \ -- \ and in fact even \  
$\bar{\vars}^{n}:=\vars^{n}\rest{\bigcup_{i=1}^{n}\Image(d^{i})}$ \ -- \ 
are epimorphisms.

Moreover, given \ 
$\alpha=\sum_{i=1}^{n} d^{k}a_{k}\in \bigcup_{i=1}^{n}\Image(d^{i})$ \ such 
that \ $s^{j}\alpha=0$ \ for \ $0\leq j\leq n-1$, \ the identities \eqref{ezero} 
imply that \ 
$$
\alpha=\sum_{p=1}^{n}\sum_{q=0}^{p-1} d^{p}d^{q}
(\sum_{i=q}^{p-2}(-1)^{p-i-1}s^{i}a_{q}+
\sum_{j=0}^{q-1}(-1)^{p-j}s^{j}(a_{q+1} + a_{p}))\ 
\in \bigcup_{i=1}^{n-1}\Image(d^{i}),
$$ 
\noindent so by induction on $n$ we see \ $\alpha=0$ \ -- \ and thus \ 
$\bar{\vars}^{n}$ \ is one-to-one, \ so in fact it is an isomorphism of 
unstable coalgebras. We thus have \ $\Aus=G(R\lra{\Xs\amalg\Ys})$ \ for some \ 
$\Xs\in\gr\Set$, \ where $\amalg$ denotes the disjoint union, and we may assume \ 
$\bigcup_{i=1}^{n}\Image(d^{i})=G(R\lra{\Ys})$.

Finally, for each \ $x\in\Xs$, \ set \ $x_{0}\DEF x$, \ and define \ $x_{k}$ \ 
inductively by \ $x_{k+1}\DEF x_{k}-d^{k+1}s^{k}x_{k}$ \ ($0\leq k<n$). \ 
We see that \ $s^{j}x_{k}=0$ \ for \ $0\leq j<k$, \ so \ $\hat{x}\DEF x_{n}$ \ 
has \ $s^{j}\hat{x}=0$ \ for all $j$, \ i.e., \ $\vars^{n}(\hat{x})=0$. \ 
Moreover, \ $x-\hat{x}\in G(R\lra{\Ys})$, \ so if we set \ 
$\bAs{n}\DEF G(R\lra{\{\hat{x}\}_{x\in\Xs}})$, \ we get the required CW basis.
\end{proof}

\subsection{The resolution model category \ $c\Sa$}
\label{srmc}\stepcounter{thm}

In \cite[\S 5.8]{DKStE}, Dwyer, Kan and Stover define a model category structure 
on the category \ $c\Sa$ \ of cosimplicial spaces (for each choice of $R$), 
which they called the ``$E^{2}$-model category'', though the term 
\emph{resolution model category} (cf.\ \cite{GHopR}) may perhaps be more 
appropriate:

A map \ $f:\Xu\to\Yu$ \ of cosimplicial spaces is

\begin{enumerate}
\renewcommand{\labelenumi}{(\roman{enumi})~}
\item a \emph{weak equivalence} if \ $\pi^{n}\HRs{f}$ \ is an isomorphism 
(of graded $R$-modules) for each \ $n\geq 0$;
\item a \emph{cofree map} if for each \ $n\geq 0$ \ there is a fibrant $R$-GEM \ 
$K^{n}\in\Sa$ \ and a map \ $X^{n}\to K^{n}$ \ which induces a trivial 
fibration \ $X^{n}\to (X^{n}\times_{L^{n}\Xu} L^{n}\Yu)\times K^{n}$;
\item a \emph{fibration} if it is a retract of a cofree map;
\item a \emph{cofibration} if \ 
$f^{n}\bot\xi^{n}:X^{n}\amalg_{M^{n}\Yu}M^{n}\Xu\to Y^{n}$ \ (cf.\ \S \ref{dmat}) 
is a cofibration for each \ $n\geq 0$, \ and \ $\pi^{n}f$ \ is a levelwise 
cofibration (i.e., monomorphism) of graded $R$-modules.
\end{enumerate}

The advantage of such a model category is that it provides a way to define a 
\emph{cosimplical resolution} of a simplicial set (or topological space) \ 
$\X\in\Sa$, \ as a fibrant replacement for the constant cosimplicial space \ 
$\co{\X}$ \ -- \ where a special case (in fact, the motivating example) is the 
$R$-resolution presented in \cite[I, \S 4.1]{BKaH}. \ See also 
\cite[\S 2]{BlaAI} for a slight generalization of the original construction.
%
%
\sect{The fundamental group and cohomology}
\label{cfc}

As noted in \S \ref{rcoh} above, the category \ $c\CA{R}$ \ of cosimplicial 
unstable $R$-coalgebras has a model category structure in which the objects \
$\SKp{n}{k}{\Ca}$ \ ($n\geq 0$, \ $k\geq 1$) \ play the role of 
cosimplicial Eilenberg-Mac Lane objects, in the sense of representing the
cohomotopy groups.  Thus, if we take homotopy classes of maps 
between (products of) these Eilenberg-Mac Lane objects as the primary cohomotopy 
operations (see \cite[V, \S 8]{GWhE}), we can endow the cohomotopy groups \ 
$\pus\Aus=(\pi^{i}\Aus)_{i=0}^{\infty}$ \ of any \ $\Aus\in c\CA{R}$ \ with an 
additional structure: that of a \emph{$\CA{R}$-coalgebra}, \ that is, a graded
object over \ $\CA{R}$ \ (coabelian, in positive dimensions), endowed with an 
action of these primary cohomotopy operations. This concept is dual to 
that of a \ $\K{R}$-$\Pi$-algebra, \ in the terminology of \cite[\S 3.2]{BStG}. \ 
By definition, this structure is a homotopy invariant of \ $\Aus$. 

\subsection{The coaction of the fundamental group}
\label{scfg}\stepcounter{thm}

In our case we shall only need the very simplest part of this structure \ -- \ 
namely, the coaction of the fundamental group \ $\pi^{0}\Aus$ \ on each of the 
higher cohomotopy groups \ $\pi^{n}\Aus$. \ This is described in terms of 
homotopy classes of maps \ 
$\SKp{n}{p}{\Ca}\to\SKp{0}{k}{\Ca}\times\SKp{n}{\ell}{\Ca}$; \ but 
since these are cosimplicial coalgebras of finite type, it is may be easier to
follow the dual description, in \ $s\K{R}$, \ of homotopy classes of maps between 
simplicial suspensions \ $\SKp{n}{k}{\K{R}}\in s\K{R}$ \ of the free
unstable algebras \ $\Kp{k}{\K{R}}\DEF\Hu{\ast}{\EM{R}{k}}{R}\in\K{R}$.

First, if \ $\Yds$ \ is any simplicial graded-commutative algebra over a field 
$R$, we can define the ``$\star$-action'' of any \ $a\in Y_{0}^{k}$ \ on \ 
$b\in Y_{n}^{\ell}$ \ by \ $a\star b\DEF \hat{a}\cdot b\in Y_{n}^{k+\ell}$, \ 
where \ $\hat{a}\DEF s_{n-1}\dotsb s_{0}a\in Y_{n}^{k}$. \  If we define 
the $n$-\emph{cycles} and $n$-\emph{chains} algebras \ 
$Z_{n}\Yds\subset C_{n}\Yds$ \ dually to \ $C^{n}\Xu\epic B^{n}\Xu$ \ of 
\S \ref{dnc} \ (see \cite[\S 17]{MayS}), \ then since $\star$ commutes with
the face maps, it defines a (bilinear) ``action'' of \ $Y_{0}^{\ast}$ \ 
on \ $C_{n}\Yds$ \ and \ $Z_{n}\Yds$ \ and thus an action of \ $\pi_{0}\Yds$ \ 
on \ $\pi_{n}\Yds$ \ for any \ $n\geq 1$.

Now let \ $\Xds$ \ denote the simplicial unstable algebra \ 
$\SKp{0}{k}{\K{R}}\times\SKp{n}{\ell}{\K{R}}\in s\K{R}$, \ for \ $k,\ell,n>0$. \ 
Note that we have a short exact sequence of unstable algebras:
\setcounter{equation}{\value{thm}}\stepcounter{subsection}
\begin{equation}\label{efour}
0\to Z_{n}\Xds\hra
\Hu{\ast}{\EM{R}{k}\times\EM{R}{\ell}}{R}\xra{\bld}
\Hu{\ast}{\EM{R}{k}}{R}\to 0.
\end{equation}
\setcounter{thm}{\value{equation}}
\noindent Evidently \ 
$Z_{n}\Xds$ \ consists of elements of the form \ $\sum_{i}a_{i}\star b_{i}$ \ 
where \ $a_{i}\in X_{0}^{\ast}$ \ and \ $0\neq b_{i}\in X_{n}^{\ast}$. \ 
However, if \ $b=b'\cdot b''$ \ is non-trivially decomposable in \ 
$X_{n}^{\ast}$, \ then \ 
$\zeta\DEF(s_{n}\hat{a})\cdot (s_{0}b'\cdot s_{0}b''-s_{0}b'\cdot s_{1}b'')
\in X_{n+1}^{\ast}$ \ 
satisfies \ $d_{0}\zeta=\hat{a}\cdot b'\cdot b''$ \ and \ $d_{j}\zeta=0$ \ for \ 
$1\leq j\leq n+1$, \ so that \ $\pi_{n}\Xds$ \ is just the free \ 
$\pi_{0}\Xds$-module generated by \ $(\Kp{\ell}{\K{R}})_{ab}$ \ (see 
\S \ref{sac}, and compare \cite[\S 5]{BStG}), \ where \ 
$\pi_{0}\Xds\cong\Kp{k}{\K{R}}=\Hu{\ast}{\EM{R}{k}}{R}\in\K{R}$. \ 

For the dual category of cosimplicial coalgebras, we need the following

\begin{defn}\label{dcm}\stepcounter{subsection}
For a given coalgebra \ $\Js\in\CA{R}$, \ an unstable coalgebra \ $\Cs$ \ 
equipped with bilinear ``co-operations'' \ $\Cs\to\Js\otimes\Cs$, \ 
(satisfying the universal identitites for the dual of ``action'' \ 
$\star:\pi_{0}\Yds\otimes\pi_{n}\Yds\to\pi_{n}\Yds$ \ defined above) will be
called a \ $\Js$-\emph{coalgebra}. \ The category of such will be denoted by \ 
$\CA{\Js}$.

On the other hand, a \textit{coabelian} unstable coalgebra \ $\Ks$ \ equipped 
with a (left) coaction map of coalgebras \ $\psi:\Ks\to\Js\otimes\Ks$, \ 
satisfying the usual identities (see \cite[\S 2.1]{SweeH}) is called a 
$\Js$-\emph{comodule}. \ We denote the category of such by \ $\CM{\Js}$. \ 
This is an abelian category.

We say that a $\Js$-comodule \ $\Ks$ \ is \emph{cofree}, with \emph{basis} \ 
$\Vs\in gr\RM$, \ if \ $\Ks=\Vs\otimes_{R}\Js$ \ (as graded $r$-modules), \ 
and the coaction \ $\psi:\Vs\otimes_{R}\Js\to(\Vs\otimes_{R}\Js)\otimes_{R}\Js$ \ 
is induced by the compultiplication \ $\Js\to\Js\otimes_{R}\Js$. \ Similarly,
a map of cofree comodules is called \emph{cofree} if it is induced by a map of
the bases.
\end{defn}

The above discussion for the case of unstable algebras may now be summarized in:
%
%
\begin{prop}\label{pthree}\stepcounter{subsection}
Any cosimplicial unstable coalgebra \ 
$\Aus\in c\CA{R}$ \ has a coaction of \ $\pi^{0}\Aus\in\CA{R}$ \ on \ 
$\pi^{n}\Aus\in(\CA{R})_{\ca}$, \ induced by an \ $\As^{0}$-coalgebra structure 
on \ $\As^{n}$. \ This coaction of \ $\As^{0}$ \ commutes with the \ 
$\MA{R}$-structure (i.e., the action of the Steenrod algebra), and respects the 
coface maps, and thus \ $\bd:B^{n}\Aus\to C^{n+1}\Aus$ \ is a map of \ 
$\As^{0}$-coalgebras for each \ $n\geq 0$.
\end{prop}

Note that the \ $\pi^{0}\Aus$-comodule structure on the coabelian coalgebra \ 
$\pi^{n}\Aus$ \ is just part of a bigraded cocommutative coalgebra structure on \ 
$\Css\DEF\pi^{\ast}\Aus$, \ in which the diagonal respects the unstable 
$R$-operations. This is the cosimplicial analogue of the \ $\C$-\Pa-structure
on the homotopy objects of a simplicial object over a category of 
universal algebras $\C$ (see \cite[\S 3.2]{BStG}).

\subsection{Quillen cohomology}
\label{sqcc}\stepcounter{thm}

Given an unstable coalgebra \ $\Js\in\CA{R}$, \ and \ $\Ks\in\CM{\Js}$ \ -- \ 
that is, a coabelian object \ $\Ks\in(\CA{R})_{\ca}$, \ with a coaction of \ 
$\Js$ \ -- \ one may define its Quillen cohomology by dualizing 
\cite[\S 2]{QuC}, as follows:

Choose some cofree cosimplicial resolution \ $\Aus\in c\CA{R}$ \ of \ $\Js$, \ 
and note that the \ $\Js$-comodule \ $\Ks$ \ is in particular an \ 
$\As^{0}$-comodule, and each \ $\As^{n}$ \ is a coalgebra over \ $\As^{0}$.
Moreover, by the usual universal property we have a natural equivalence \ 
$\Hom_{\CA{\As^{0}}}(\Ks,\Aus)\cong\Hom_{\CM{\As^{0}}}(\Ks,(\Aus)_{\ca'})$, \ 
where \ $M_{\ca'}$ \ denotes the coabelianization of the $\As^{0}$-coalgebra 
$M$ \ (see \S \ref{sac}). \ Thus \ $\Cud\DEF\Hom_{\CA{\As^{0}}}(\Ks,\Aus)$ \ is a 
cosimplicial abelian group, and \ $\pi^{n}\Cud$ \ is called the $n$-th 
\emph{Quillen cohomology group} of \ $\Js$ \ \emph{with coefficients in} \ 
$\Ks$, \ and denoted by \ $\Hu{n}{\Js}{\Ks}$ \ (compare \cite[\S 3]{KLeeC}). \ 
If \ $\Js$, \ $\Ks$ \ and \ $\Aus$ \ are of finite type, this is the vector 
space dual of the usual Quillen cohomology of \ $(\Js)^{\star}\in\K{R}$ \ 
(see \cite[\S 8]{BousH}, and compare \cite[\S 4]{BlaAI}).

\begin{remark}\label{rqc}\stepcounter{subsection}
As for any abelian category, given a cosimplicial object \ 
$(\Aus)_{\ca'}$ \ over \ $\CM{\As^{0}}$, \ in each dimension \ $n\geq 1$ \ there
is a direct product decomposition \ 
$(\As^{n})_{\ca'}=C^{n}((\Aus)_{\ca'})\oplus L^{n}((\Aus)_{\ca'})$ \ (compare \ 
\cite[Cor.\ 22.2]{MayS}). If we choose a CW basis \ $(\bAs{n})_{n=0}^{\infty}$ \ 
for \ $\Aus$ \ (\S \ref{dcw}), \ we have:
$C^{n}(\Aus)_{\ca'}\cong(\bAs{n})_{\ca}\otimes\As^{0}$ \ 
as unstable modules (where \ $(\bAs{n})_{\ca}$ \ is the usual coabelianization 
of \S \ref{sac}).  Moreover, each \ $(\As^{n})_{\ca'}$ \ is a cofree \ 
$\As^{0}$-comodule (with a basis which may be described explicitly in terms
of the $CW$-basis for \ $\Aus$ \ -- \ see \cite[(6.3)]{BlaN}), and the coface maps
are cofree (\S \ref{dcm}), so we have \ 
$C^{n}((\Aus)_{\ca'})=(C^{n}\Aus)_{\ca'}$. \ We may therefore use the cochain 
complex \ 
$$
\dotsb\to \Hom_{\CM{\As^{0}}}(\Ks,(C^{n}\Aus)_{\ca'})\xra{\delta^{n}}
\Hom_{\CM{\As^{0}}}(\Ks,(C^{n+1}\Aus)_{\ca'})\to\dotsb
$$
\noindent (where \ $\delta^{n}$ \ is induced by \ 
$\dd{A^{n}}\circ q^{n}:C^{n}\Aus\to C^{n+1}\Aus$) \ to calculate \ 
$\Hu{\ast}{\Js}{\Ks}$ \ (compare \cite[X, \S 7.1]{BKaH}).
\end{remark}

\subsection{An alternative description}
\label{sad}\stepcounter{thm}

Quillen's original description of the cohomology of a (simplicial) algebra
included several variant approaches (cf.\ \cite[\S 3]{QuC}, and by dualizing
one of these, Bousfield obtained an alternative description of the cohomology of a
coalgebra, as follows:

For any field $R$, given an unstable coalgebra \ $\Js\in\CA{R}$ \ and a
$\Js$-comodule \ $\Ks\in\CM{\Js}$, \ with coaction \ 
$\psi:\Ks\to\Js\otimes\Ks$, \ and a map of coalgebras \ $\delta:\Js\to\Ls$, \ 
define a \emph{derivation} \ $\varphi:\Ks\to\Ls$ \ to be a \ 
$\MA{R}$-morphism such that \ 
$\Delta_{\Js}\circ\varphi=
\delta\otimes \varphi+\tau\circ(\varphi\otimes\delta)\circ\psi$ \ 
(see \S \ref{dca}). \ Write \ 
$\Der_{\CA{R}}(\Ks,\Ls)$ \ for the $R$-vector space of all such derivations.

For every comodule \ $\Ks\in\CM{\Js}$, \ we can think of \ $\Js\oplus\Ks$ \ 
as a coalgebra under \ $\Js$, \ with diagonal \ $\Delta_{\Js\oplus\Ks}$ \ 
defined by\ 
$$
\Delta_{\Js}+\psi+\tau\circ\psi+0:\Js\oplus\Ks\to
(\Js\otimes\Js)\oplus(\Js\otimes\Ks)\oplus(\Ks\otimes\Js)\oplus(\Ks\otimes\Ks).
$$ \ 
Thus, given a map \ $\delta:\Js\to\Ls$, \ we have a natural identification
$$
\Hom_{\CA{R}\backslash\Js}(\Js\oplus\Ks,\Ls)\cong \Der_{\CA{R}}(\Ks,\Ls).
$$

In fact, the functor \ $\Ks\mapsto\Js\oplus\Ks$ \ is left adjoint to the 
coabelianization functor \ $(\,)_{\ca}:\CA{R}\backslash\Js\to\CM{\Js}$ \ 
(compare \S \ref{sac}), and it induces an equivalence of categories between \ 
$\CM{\Js}$ \ and \ $(\CA{R}\backslash\Js)_{\ca}$. \ 
Moreover, as in \cite[\S 4]{QuC}, one has an explicit description of the 
functor \ $(\,)_{\ca}$ \ in terms of a cotensor product of \ $\Js$ \ with \ 
$\Omega\Ls$ \ (the coalgebraic analogue of the usual K\"{a}hler module of 
differentials).  
See \cite[\S 8.5]{BousH} and \cite[\S 7.7-8]{SchwU} for more details on this 
approach.

Now if \ $\Js\to\Aus$ \ is a resolution, we can show that there is a natural map 
$$
\Hom_{\CA{\As^{0}}}(\Ks,C^{n}\Aus)\to
\Hom_{\CA{R}\backslash\Js}(\Js\oplus\Ks,\As^{n}),
$$
\noindent which induces an isomorphism in cohomology, so Remark \ref{rqc} above 
implies that the Quillen cohomology groups \ $\Hu{n}{\Js}{\Ks}$ \ coincide with 
the derived functors of \ $\Der_{\CA{R}}(\Ks,-)$ \ applied to \ $\Js$, \ which
were considered by Bousfield (who showed in \cite[\S 9]{BousH} that the groups \ 
$\Der^{s,t}(\Js,\Ls)\DEF\Hu{s}{\Js}{\tHi{\ast}{\bS{t}}{\Fp}\otimes\Ls}$ \ 
serve as the \ $E_{2}$-term of a certain unstable Adams spectral sequence).

To show this, use the fact that a morphism in either \ $\Hom$-set 
must take values in a coabelian unstable coalgebra, and that the unique
iterated coface map \ $\delta:\Js\to\As^{n}$ \ vanishes when projected to \ 
$C^{n}\Aus$. \ We omit the details, since we shall not require this result in 
what follows. However, it may be observed that in the Massey-Peterson case \ 
$\Js=U(\Ms)$, \ Bousfield's approach allows us to identify \ 
$\Hu{n}{\Js}{\Ks}$ \ with the usual \ $\Ext^{n}_{\MA{\Fp}}(\Ks,\Ms,)$, \ 
so we can recover Harper's results in \cite[Prop.\ 4.2]{HarpC1} as a 
special case of Theorem \ref{ttwo} below.

It is possible that one could prove the results of the following sections, using 
Bousfield's description of cohomology as the derived functors of derivations,
and dualizing Quillen's identification of these derived functors (for
algebras) with suitable groups of extensions (see \cite[\S 3]{QuC}). However,
there may be computational advantages to the explicit approach we have taken here.

%
%
\sect{Realizing resolutions}
\label{crr}

The key to our approach to the realization question for an unstable coalgebra \ 
$\Ks$ \  lies in the realization of a suitable cofree resolution of \ 
$\Ks$ \ -- \ by analogy with the method used in \cite{BlaAI} for $\Pi$-algebras.
In what follows \ $R=\Q$ \ or \ $\Fp$.

\subsection{Trying to realize a resolution}
\label{strr}\stepcounter{thm}

Given an unstable coalgebra \ $\Ks\in\CA{R}$, \ choose some cosimplicial 
resolution \ $\Ks\to\Aus$, \ with CW basis \ $\{\bAs{n}\}_{n=0}^{\infty}$, \ 
as in \S \ref{sccr}. We would like to realize this algebraic resolution at 
the space level, i.e., find a cofibrant cosimplicial space \ $\Qu\in c\Sa$, \ 
with a CW basis \ $\{\baQ{n}\}_{n=0}^{\infty}$, \ such that \ 
$\HRs{\baQ{n}}\cong\bAs{n}$ \ for all \ $n\geq 0$, \ and the attaching maps \ 
$\bar{d}_{Q^{n}}^{0}:\baQ{n}\to\baQ{n+1}$ \ realize those for \ $\Aus$, \ so 
that \ $\HRs{\Qu}\cong\Aus$. \ 

We attempt to construct such a \ $\Qu$ \ (with its CW basis) by induction on 
the cosimplicial dimension:

The first two steps are always possible: because \ $\bAs{0}\in\CA{R}$ \ 
is cofree by assumption (\S \ref{dcw}), we can find a map \ 
$\bar{d}_{Q}^{0}:\baQ{0}\to\baQ{1}$ \ in \ $\Sa$ \ 
such that \ $(\bar{d}_{Q}^{0})_{\#}:\HRs{\baQ{0}}\to\HRs{\baQ{1}}$ \ is \ 
$\bar{d}^{0}_{A}:\bAs{0}\to\bAs{1}$ \ (\S \ref{rcuc}).
We then set \ $\bQ^{0}\DEF\baQ{0}$ \ and \ 
$\tilde{\bQ}^{1}\DEF\baQ{1}\times L^{1}\Qu=\baQ{1}\times\bQ^{0}$, \ with \ 
$\tilde{d}^{0}\DEF\bar{d}^{0}\top\Id$ \ and \ $\tilde{d}^{1}\DEF0\top\Id$. \ 
In order to end up with a cofibrant cosimplicial space (\S \ref{dmat}), \ we now 
change the resulting \ $\tilde{\xi}^{n}:M^{1}\Qu\to\tilde{\bQ}^{1}$ \ into a 
cofibration: \ $\xi^{n}:M^{1}\Qu\to\bQ^{1}$. \ It will be convenient to
denote \ $d^{0}\circ p^{1}:\bQ^{0}\to C^{1}\Qu$ \ by \ 
$\bd=\dd{Q^{0}}:B^{0}\Qu\to C^{1}\Qu$, \ to conform with the notation of Figure
\ref{fig1}.

If we let \ $B^{1}\Qu$ \ denote the cofiber of \ $\dd{Q^{0}}$, \ from the long 
exact sequence in homology:
$$
\dotsc  \HR{i+1}{B^{1}\Qu}\xra{\partial_{1}}\HR{i}{B^{0}\Qu}
\xra{\bd}\HR{i}{C^{1}\Qu}\xra{\qu{1}{Q}}\HR{i}{B^{1}\Qu}\dotsc
$$
\noindent and Proposition \ref{pone} we obtain a short exact sequence of 
unstable coalgebras:
$$
0\to B^{1}(\Aus)\xra{i}\HRs{B^{1}\Qu}\xepic{\psi}\Sigma\Ks\to 0,
$$

\noindent where \ $\Sigma\Ks$ \ is the graded $R$-module, shifted one
degree up, with the trivial coal\-ge\-bra structure (and the suspended action
of the Steenrod algebra, if \ $R=\Fp$)\vsm.

Now assume \ $\Qu$ \ as required has been constructed through cosimplicial 
dimension $n$, \ so we have an $n$-cosimplicial object which (by a slight abuse 
of notation) we denote by \ 
$\tr{n}\Qu\in c^{\lra{n}}\Sa$ \ (cf.\ \S \ref{dcs}), with \ 
$\HRs{\tr{n}\Qu}\cong\tr{n}\Aus$. \ 

By considering the cosimplicial chain complex corresponding to \ $R\lra{\Qu}$, \ 
one can verify (as in the proof of Proposition \ref{pone}) that there always 
exists a factorization of \ $\qu{n}{Q}$ \ as follows:

%
%
\begin{figure}[htbp]
\begin{picture}(180,60)(-70,0)
%
%
\put(20,50){$C^{n}\Aus$}
\put(55,55){\vector(1,0){75}}
\put(130,55){\vector(1,0){3}}
\put(85,62){$q^{n}_{A}$}
\put(140,50){$B^{n}\Aus$}
\put(175,55){\vector(1,0){3}}
\put(178,55){\vector(1,0){75}}
\put(207,62){$\dd{A^{n}}$}
\put(260,50){$C^{n+1}\Aus$}
%
%
\put(35,47){\vector(0,-1){33}}
\put(20,30){$\cong$}
\multiput(155,47)(0,-3){10}{\circle*{.3}}
\put(155,15){\vector(0,-1){2}}
\put(158,30){$\bar{q}^{n}$}
\put(280,47){\vector(0,-1){33}}
\put(285,30){$\cong$}
%
%
\put(0,0){$\HRs{C^{n}\Qu}$}
\put(70,5){\vector(1,0){45}}
\put(80,11){$\qu{n}{Q}$}
\put(120,0){$\HRs{B^{n}\Qu}$}
\put(195,5){\vector(1,0){50}}
\put(207,11){$\ddi{Q^{n}}$}
\put(250,0){$\HRs{C^{n+1}\Qu}$,}
\end{picture}
\caption[fig2]{}
\label{fig2}
\end{figure}

\noindent and since \ $\Aus$ \ is a resolution, \ $\dd{A^{n}}$, \ and thus \ 
$\bar{q}^{n}$, \ must be monic. Moreover, for any \ 
$a\in\HR{i}{B^{n-1}\Qu}$ \ we have \ 
$\qu{n}{Q}(\ddi{Q^{n-1}}(a))=0$; \ thus \ 
$\ddi{Q^{n-1}}(a)\in \HR{i}{C^{n}\Qu}=C^{n}\Aus$ \ is a cocyle for 
the resolution \ $\Aus$, \ so it must be in \ $\Image(\dd{A^{n-1}})$. \ Thus
we must have
\setcounter{equation}{\value{thm}}\stepcounter{subsection}
\begin{equation}\label{efive}
\Image\qu{n}{Q} + \Ker\ddi{Q^{n}} = \HRs{B^{n}\Qu}.
\end{equation}
\setcounter{thm}{\value{equation}}
\noindent 

We therefore assume by induction that we have a short exact sequence of unstable 
coalgebras: \ 
\setcounter{equation}{\value{thm}}\stepcounter{subsection}
\begin{equation}\label{esix}
0\to B^{n}\Aus\xra{\bar{q}^{n}}\HRs{B^{n}\Qu}\xepic{\psi^{n}}
\Coker(\bar{q}^{n})\to 0,
\end{equation}
\setcounter{thm}{\value{equation}}
\noindent where \ $\Coker(\bar{q}^{n})\cong\Sigma^{n}\Ks$, \ and in fact \ 
$\psi^{n}\rest{\Image(\partial_{n+1})}$ \ is an 
isomorphism onto \ $\Coker(\bar{q}^{n})$, \ with \ $\partial_{n+1}$ \ 
again the connecting homomorphism in the long exact sequence
\setcounter{equation}{\value{thm}}\stepcounter{subsection}
\begin{equation}\label{eseven}
\dotsc  H_{i+1} B^{n+1}\Qu \xra{\partial_{n+1}} H_{i} B^{n}\Qu
\xra{\bd} H_{i} C^{n}\Qu \xra{(q^{1})_{\#}} H_{i} B^{n}\Qu
\xra{\partial_{n+1}}\dotsc 
\end{equation}
\setcounter{thm}{\value{equation}}

Since \ $\Image\qu{n}{Q}=\Image\bar{q}^{n}$ \ in
Figure \ref{fig2}, in fact we have a direct sum of $R$-modules in \eqref{efive}, 
and this is by definition a semi-split extension of coalgebras. This implies
that
\setcounter{equation}{\value{thm}}\stepcounter{subsection}
\begin{equation}\label{eeight}
\partial_{n-1}\rest{\Image(\partial_{n}) } \ \ \ 
\text{is one-to-one, and surjects onto} \ \Image(\partial_{n-1}).
\end{equation}
\setcounter{thm}{\value{equation}}

\begin{aside}\label{apost}\stepcounter{subsection}
Observe that if we contine our construction ``naively''
by choosing some GEM \ $\baQ{n+1}\in\Sa$ \ with an attaching map \ 
$\bar{d}^{0}_{n+1}:\bQ^{n+1}\to\baQ{n+1}$ \ which induces a monomorphism \ 
$\HRs{B^{n}\Qu}\hra \HRs{\baQ{n+1}}$ \ in homology, we can easily continue
this process to obtain a cosimplicial space \ $\Yu$, \ such that \ 
$\HRs{\Yu}$ \ is a free cosimplicial coalgebra satisfying:
\setcounter{equation}{\value{thm}}\stepcounter{subsection}
\begin{equation}\label{enine}
\pi^{i}\HRs{\Yu}\cong 
\begin{cases}\Ks & \text{for \ }i=0,\\
\Sigma^{n}\Ks    & \text{for \ }i=n+1,\\
0                & \text{otherwise}. \end{cases}
\end{equation}
\setcounter{thm}{\value{equation}}

Such a \ $\Yu$ \ should be thought of as the $(n-1)$-\emph{st Postnikov section} 
for the resolution \ $\Qu$. \ We denote a cofibrant version of it by \ 
$\bP{n-1}\Qu$, \ and observe that it is unique up to homotopy equivalence 
(in the model category structure of \S \ref{srmc}). 
This provides a convenient homotopy-invariant version of the $n$-coskeleton of
a cosimplicial space. See \cite[\S 3.4]{BlaAI} for an explanation of the indexing.

In \cite{BGoeC}, we present an alternative approach to the dual problem of
realizing (simplicial) \Pa s, via Postnikov systems (including 
$k$-invariants) for simplicial \Pa s and simplicial spaces. This was in fact the 
original program of \cite{DKStE,DKStB}. However, because 
there is no satisfactory homotopy theory for cosimplicial \emph{sets}, an 
analogous approach here would require developing additional machinery not 
presently available.
\end{aside}

\subsection{Continuing the construction}
\label{scc}\stepcounter{thm}

If we can extend our $n$-truncated object \ $\tr{n}\Qu$ \ one more dimension, 
we will have a principal face map \ $\dd{n}:B^{n}\Qu\to C^{n+1}\Qu$, \ which 
induces a map \ 
$\lambda=\ddi{n}:\HRs{B^{n}\Qu}\to\HRs{C^{n+1}\Qu}\cong C^{n+1}\Aus$ \ (by
Proposition \ref{pone}). It turns out that such a $\lambda$ is essentially
all we need in order to proceed.

First note that since \ $\bAs{n+1}$ \ is a cofree unstable coalgebra by 
assumption, and \ $\bar{q}^{n}:B^{n}\Aus\to\HRs{B^{n}\Qu}$ \ 
in Figure \ref{fig2} above is monic, \ $\bbd_{A^{n}}:B^{n}\Aus\to\bAs{n+1}$ \ 
extends to a map \ $\lambda:\HRs{B^{n}\Qu}\to\bAs{n+1}$ \ as follows:

%
%
\begin{figure}[htbp]
\begin{picture}(180,120)(-100,5)
%
%
\put(0,100){$B^{n}\HRs{\Qu}$}
\put(72,105){\vector(1,0){2}}
\put(74,105){\vector(1,0){42}}
\put(85,110){$\bar{q}^{n}$}
\put(120,100){$\HRs{B^{n}\Qu}$}
\put(195,105){\vector(1,0){40}}
\put(235,105){\vector(1,0){2}}
\put(205,109){$\psi^{n}$}
\put(240,100){$\Sigma^{n}\Ks$}
%
%
\put(32,94){\vector(0,-1){28}}
\put(20,75){$\cong$}
\multiput(170,93)(3,-3){8}{\circle*{.3}}
\put(194,69){\vector(1,-1){2}}
\put(190,83){$\lambda$}
%
%
\put(20,55){$B^{n}\Aus$}
\put(55,60){\vector(1,0){3}}
\put(58,60){\vector(1,0){45}}
\put(70,65){$\dd{A}$}
\put(105,55){$C^{n+1}\Aus$}
\put(148,60){\vector(1,0){42}}
\put(190,60){\vector(1,0){3}}
\put(155,65){$\proj$}
\put(195,55){$\bAs{n+1}$}
%
%
\put(25,5){$\As^{n}$}
\put(40,13){\vector(4,1){157}}
\put(110,17){$\bar{d}^{0}_{A}$}
\put(32,18){\vector(0,1){31}}
\put(32,49){\vector(0,1){3}}
\put(-5,30){$q^{n}_{A}\circ p^{n}_{A}$}
\end{picture}
\caption[fig3]{}
\label{fig3}
\end{figure}

We can realize $\lambda$ by a map \ $\ell:B^{n}\Qu\to\baQ{n+1}$, \ 
for a suitable GEM \ $\baQ{n+1}$, \ and thus a map \ 
$\bar{d}^{0}:\bQ^{n}\to\baQ{n+1}$ \ (see Figure \ref{fig1}), \ which in turn 
determines an extension of \ $\tr{n}\Qu$ \ to \ $\tr{n+1}\Qu$; \ we may
modify this to be cofibrant. Moreover, from the exactness of the bottom row,
by \eqref{efive}, we have a unique lifting \ 
$\mu:\Coker(\bar{q}^{n})\to B^{n+1}\Aus$ \ as follows:

%
%
\begin{figure}[htbp]
\begin{picture}(340,120)(-50,-5)
%
%
\put(0,100){$B^{n}\HRs{\Qu}$}
\put(72,105){\vector(1,0){2}}
\put(74,105){\vector(1,0){32}}
\put(85,111){$\bar{q}^{n}$}
\put(110,100){$\HRs{B^{n}\Qu}$}
\put(185,105){\vector(1,0){40}}
\put(220,105){\vector(1,0){2}}
\put(195,109){$\psi^{n}$}
\put(230,100){$\Coker(\bar{q}^{n})\cong\Sigma^{n}\Ks$}
%
%
\put(32,95){\vector(0,-1){82}}
\put(20,55){$\cong$}
\put(145,95){\vector(0,-1){45}}
\put(148,73){$\ddi{Q}$}
\put(120,40){$\HRs{C^{n+1}\Qu}$}
\put(145,35){\vector(0,-1){20}}
\put(148,23){$\cong$}
\multiput(245,95)(0,-3){26}{\circle*{.3}}
\put(245,16){\vector(0,-1){2}}
\put(249,55){$\mu$}
\multiput(315,95)(1,-3){26}{\circle*{.3}}
\put(340,20){\vector(1,-3){2}}
\put(332,55){$\xi$}
%
%
%
\put(20,0){$B^{n}\Aus$}
\put(55,5){\vector(1,0){3}}
\put(58,5){\vector(1,0){65}}
\put(75,12){$\dd{A^{n}}$}
\put(125,0){$C^{n+1}\Aus$}
\put(168,5){\vector(1,0){62}}
\put(230,5){\vector(1,0){3}}
\put(185,12){$q^{n+1}_{A}$}
\put(235,0){$B^{n+1}\Aus$}
\put(280,5){\vector(1,0){37}}
\put(286,12){$\dd{A^{n+1}}$}
\put(320,0){$C^{n+2}\Aus$}
\end{picture}
\caption[fig4]{}
\label{fig4}
\end{figure}

\noindent and if \ $\mu=0$, \ then \eqref{efive} will hold for \ $n+1$.

\begin{defn}\label{dcc}\stepcounter{subsection}
The cohomology class \ $\chi_{n}\in H^{n+2}(\Ks;\Sigma^{n}\Ks)$ \ represented
by the cocycle \ 
$\xi\DEF\dd{A^{n+1}}\circ\mu\in\Hom_{\CA{\As^{0}}}(\Sigma^{n}\Ks,\Aus)$ \ 
(see \S \ref{rqc}) \ is called the \emph{characteristic class of the extension} 
\eqref{esix} (compare \cite[IV, \S 5]{MacH} and \cite[\S 4.5]{BlaAI}).
\end{defn}

%
%
\begin{thm}\label{tone}\stepcounter{subsection}
The cohomology class \ $\chi_{n}\in H^{n+2}(\Ks;\Sigma^{n}\Ks)$ \ is independent 
of the choice of lifting $\lambda$, \ and \ $\chi_{n}=0$ \ if and only if one 
can extend \ $\bP{n-1}\Qu$ \ to an \ $n$-th Postnikov approximation \ 
$\bP{n}\Qu$ \ of a resolution of \ $\Ks$.
\end{thm}

\begin{proof}
Assume that we want to replace $\lambda$ by a different lifting \ 
$\lambda':\HRs{B^{n}\Qu}\to\bAs{n+1}$, \ and choose maps \ 
$\ell,\ell':B^{n}\Qu\to\baQ{n+1}$ \ realizing \ $\lambda$, \ 
$\lambda'$ \ respectively; \ their respective extensions to  \ 
$d^{0}$ \ and \ $(d^{0})':\bQ^{n}\to\bQ^{n+1}$ \ agree on \ $L^{n+1}\Qu$.
We correspondingly have \ 
$\mu':\Sigma^{n}\Ks\to B^{n+1}\Aus$ \ and \ $\xi'\DEF\dd{A^{n+1}}\circ\mu'$ \ 
in Figure \ref{fig4}.

Since \ $\bQ^{n+1}$ \ can be any fibrant GEM realizing \ $A^{n+1}$, \ we may 
assume it is a simplicial $R$-module, and thus \ 
$\Hom_{\Sa}(\bQ^{n},\bQ^{n+1})$ \ has a natural $R$-module structure. Set \ 
$h\DEF ((d^{0})'-d^{0}):\bQ^{n}\to\bQ^{n+1}$. \ 
Then $h$ induces a map \ $\eta:\HRs{B^{n}\Qu}\to C^{n+1}\Aus$ \ whose projection
onto \ $\bAs{n+1}$ \ is \ $\lambda'-\lambda$. \ Moreover, because \ 
$d^{0}$ \ and \ $(d^{0})'$ \ agree with \ $\dd{A^{n}}$ \ when pulled back to \ 
$\HRs{C^{n}\Qu}$, \ we have \ $\eta\circ \bar{q}^{n}=0$, \ and thus $\eta$ 
factors through \ $\zeta:\Sigma^{n}\Ks\to C^{n+1}\Aus$, \ and this is a map of \ 
$\As^{0}$-coalgebras because \ $\Sigma^{n}\Ks$ \ is a coabelian \ 
$\As^{0}$-comodule (actually, a \ $\Ks$-comodule), and $\zeta$ is induced by 
group operations from the \ $\As^{0}$-coalgebra maps \ $\bd$ \ and \ $(\bd)'$. \ 
See Figure \ref{fig5} below.

Moreover, in the abelian group structure on \ 
$\Hom_{\CM{\Ks}}(\Sigma^{n}\Ks,-)$ \ we have \ 
$\xi'-\xi=\dd{A^{n+1}}\circ(\mu'-\mu)=\delta^{n+1}(\zeta)$ \ 
(see \S \ref{rqc}), \ so this is a coboundary, which proves independence of the 
choice of $\lambda$.

Now assume that there exists \ $\Yu\simeq \bP{n}\Qu$ \ (\S \ref{apost}) with \ 
$\tr{n}\Yu\cong\tr{n}\Qu$. \ By the discussion in \S \ref{strr}, we know that \ 
\eqref{efive} is a direct sum  for \ $n+1$, \ and since \ $\bAs{n+1}$ \ is
cofree, we can choose \ $\lambda:\HRs{B^{n}\Qu}\to\bAs{n+1}$ \ in Figure
\ref{fig3} to extend \ $\dd{A}$ \ by zero, so \ $\mu=0$ \ in Figure \ref{fig4},
and thus \ $\xi=0$.

Conversely, if \ $\chi_{n}=0$, \ we can represent it by a coboundary \ 
$\xi=\dd{A^{n+2}}\circ\vartheta$ \ for some \ $\As^{0}$-coalgebra map \ 
$\vartheta:\Ks\to C^{n+1}\Aus$, \ and thus get \ 
$\proj_{\bAs{n+1}}\circ\vartheta\circ\psi^{n}:\HRs{B^{n}\Qu}\to\bAs{n+1}$, \ 
for \ $\proj_{\bAs{n+1}}:\As^{n+1}\to\bAs{n+1}$ \ the projection. \ If we set \ 
$\lambda'\DEF \lambda-\proj_{\bAs{n+1}}\circ\vartheta\circ\psi^{n}$ \ (we can
subtract maps, because \ $\bAs{n+1}$ \ is a graded $R$-module), we have \ 
$\Image\qu{n}{Q} + \Ker\lambda' = \HRs{B^{n}\Qu}$. \ We can therefore choose \ 
$\bar{d}^{0}_{Q^{n+1}}:B^{n}\Qu\to\baQ{n+1}$ \ realizing \ $\lambda'$, \ 
and then \ $\mu'=0$, \ so that \ $\tr{n+1}\Qu$ \ 
so constructed yields \ $\bP{n}\Qu$, \  as required.
\end{proof}

\begin{notation}\label{ndl}\stepcounter{subsection}
If we wish to emphasize the dependence on the choice of $\lambda$, \ we shall
write \ $\bP{n}\Qu[\lambda]$ \ for the extension of \ $\bP{n-1}\Qu$ \ so 
constructed, \ and write \ 
$\chi_{n+1}(\lambda)\in H^{n+3}(\Ks;\Sigma^{n+1}\Ks)$ \  for the next
cohomology class (which \emph{does} depend on $\lambda$, in principle). \ 
\end{notation}

\begin{remark}\label{rdl}\stepcounter{subsection}
Note that if \ $\chi_{n}=0$, \ the choice of $\lambda$ determines the \ 
$\As^{0}$-comodule structure on \ $\Sigma^{n+1}\Ks$ \ via \eqref{esix} \ for \ 
$n+1$. \ Moreover, for each \ $n\geq 1$, the resulting coaction \ 
$\psi_{n}:\Sigma^{n}\Ks\to\Ks\otimes\Sigma^{n}\Ks$ \ in fact agrees with the 
obvious \ $\Ks$-comodule structure, defined via the original comultiplication \ 
$\Delta:\Ks\to\Ks\otimes\Ks$: \ that is, if \ 
$\Delta(a)=\sum_{i}a'_{i}\otimes a''_{i}$, \ and \ 
$\sigma_{n}:\Ks\to\Sigma^{n}\Ks$ \ is the re-indexing isomorphism (in \ 
$\MA{R}$), \ then \ 
$\psi_{n}(\sigma_{n}(a))=\sum_{i}a'_{i}\otimes\sigma_{n}(a''_{i})$. \ This follows
from the description in \S \ref{scfg}, and the fact that the exact sequences
\eqref{eseven} (and thus also \eqref{esix}) respect the $\As^{0}$-coalgebra 
structure (Proposition \ref{pthree}).
\end{remark}

\begin{defn}\label{dcv}\stepcounter{subsection}
Note also that by standard homotopical algebra arguments the elements \ 
$\chi_{n}\in H^{n+2}(\Ks;\Sigma^{n}\Ks)$ \ do not depend on the choice of 
resolution \ $\Ks\to\Aus$. \ 
If for some (and thus any) cosimplicial cofree resolution \ $\Ks\to\Aus$, \ 
there are successive choices of liftings \ $(\lambda_{n})_{n=0}^{\infty}$ \ 
in Figure \ref{fig3} such that \ $\chi_{n+1}(\lambda_{n})=0$ \ for \ 
$n=0,1,\dotsc$, \ we say that we have a \emph{coherently vanishing} sequence of 
characteristic classes. 
\end{defn}

Thus we may encapsulate our results so far in 
%
%
\begin{cor}\label{cone}\stepcounter{subsection}
Any cofree cosimplicial resolution \ $\Aus$ \ of an unstable coalgebra \ 
$\Ks\in\CA{R}$ \ is realizable by a cosimplicial space \ $\Qu\in c\Sa$ \ 
(with \ $\Aus\cong\HRs{\Qu}$) \ if and only if \ $\Ks$ \ has a coherently 
vanishing sequence of characteristic classes.
\end{cor}

%
%
\sect{Realizing coalgebras}
\label{crec}

We now apply Theorem \ref{tone}, on the realization of cosimplicial 
resolutions of coalgebras, to the original question, namely, that of realizing
a given abstract coalgebra as the cohomology of a space.   It turns out that
the obstructions described in the previous section are all that is needed, at 
least in the simply-connected case.

\subsection{The homology spectral sequence}
\label{shss}\stepcounter{thm}

In \cite[\S 3]{RecS} and \cite{AndG}, Rector and Anderson defined the 
homology spectral sequence of a cosimplicial space \ $\Yu$ \ (see also
\cite[\S 2]{BousHS}). \ This is a second quadrant spectral sequence with \ 
\setcounter{equation}{\value{thm}}\stepcounter{subsection}
\begin{equation}\label{eten}
E^{2}_{p,q}\cong \pi^{p}\HR{q}{\Yu}
\end{equation}
\setcounter{thm}{\value{equation}}
\noindent abutting to \ $\HRs{\Tot\Yu}$, \ where the \emph{total space} \ 
$\Tot\Yu\in\Ss$ \ of a cosimplicial space \ $\Yu\in c\Ss$ \ is defined 
(cf.\ \cite[I, \S 3]{BKaH}) to be the simplicial set \ $\Td\in\Ss$ \ with \ 
$T_{n}=\Hom_{c\Ss}(\Delta[q]\times\Du, \Yu)$ \ (see \S \ref{ecs}).

In general, this spectral sequence need not converge. However, under rather 
special conditions one does have strong convergence (see \cite{BousHS} and 
\cite{ShipC}), and this yields the following:
%
%
\begin{thm}\label{ttwo}\stepcounter{subsection}
For \ $R=\Q$ \ or \ $\Fp$, \ a simply-connected unstable 
$R$-coalgebra \ $\Ks$ \ is realizable as the homology of some simply-connected
space \ $\X\in\Ta$ \ if and only if \ $\Ks$ \ has a coherently vanishing sequence 
of characteristic classes.
\end{thm}

\begin{proof}
First note that any simply-connected coalgebra over $\Q$ is realizable by 
\cite[Thm.\ I]{QuR} and its Corollary, so in this case the theorem merely states 
that one always has a coherently vanishing sequence of characteristic classes.

Given any connected space \ $\X\in\Sa$, \ the cosimplicial space \ 
$\Yu$ \ defined by \ $\Y^{n}=\bar{R}^{n+1}\,\X$ \ (where \ $\bar{R}:\Sa\to\Sa$ \ 
is the Bousfield-Kan monad \ 
$(\bar{R}\,\X)_{k}\DEF\{\sum_{i}r_{i}x_{i}\in R\lra{X_{k}}~|\ 
\sum_{i}r_{i}=1\}$ \ 
-- \ cf.\ \cite[I, \S 2]{BKaH}), \ is a cosimplicial resolution of \ $\co{\X}$ \ 
in the sense of \S \ref{srmc} \ -- \ i.e., \ $\HRs\Yu$ \ is a cosimplicial
cofree resolution of \ $\HRs{\X}$ \ (see \cite[11.5]{BKaS}). \ But then
by Corollary \ref{cone}, \ $\HRs{\X}$ \ has a coherently vanishing sequence 
of characteristic classes. 

Conversely, assume that \ $\Ks$ \ is a simply-connected unstable $R$-coalgebra
with a coherently vanishing sequence of characteristic classes.
By Corollary \ref{cone}, \emph{any} cosimplicial cofree resolution \ 
$\Ks\to\Aus$ \ may be realized by a cosimplicial space \ $\Qu\in c\Sa$. \ 
In particular, since \ $K_{0}=R$ \ and  \ $K_{1}=0$, \ we may assume that 
the same holds for each \ $\As^{n}$, \ so that each $R$-GEM \ $\bQ^{n}$ \ is 
simply-connected. \ Because \ $\Aus$ \ is a resolution, \ 
$\pi^{n}\tHi{n+s}{\Qu}{R}=0$ \ for \ $n=0$ and \ $s\leq 1$ \ or \ 
$n\geq 1$, \ and thus by \cite[Thm.\ 3.4]{BousHS} the homology spectral sequence
for \ $\Qu$ \ converges strongly to \ $\HRs{\Tot\Qu}$. \ Since the \ $E^{2}$-term
of \eqref{eten} is concentrated along the $0$-line, we get \ 
$$
\Ks\cong \pi^{0}\Aus=E^{2}_{0,\ast}\xra{\cong}\HRs{\Tot\Qu}, 
$$
\noindent and this is an isomorphism of unstable coalgebras, since the edge 
homomorphism is induced by a topological map \ $\Tot\Qu\to\Tot_{0}\Qu=\bQ^{0}$.
\end{proof}

\subsection{The non simply-connected case}
\label{snsc}\stepcounter{thm}

The simple-connectivity of \ $\Ks$ \ was only needed to guarantee convergence of
the homology spectral sequence, using \cite[Thm.\ 3.4]{BousHS}; the (algebraic) 
obstruction theory described in the previous section is of course also valid 
in the non simply-connected case.
In particular, by dualizing \cite[Prop.\ 5.1.4]{BlaD} we can construct a 
resolution with CW basis \ $(\bAs{n})_{n=1}^{\infty}$ \ with strictly increasing 
connectivity \ -- \ so that in particular \ $\bar{A}_{s}^{n}=0$ \ for \ 
$0<s\leq n$ \ -- \ and then realize \ $\Aus$ \ by a cosimplicial space \ $\Qu$, \ 
assuming that \ $\Ks$ \ has a coherently vanishing sequence of characteristic 
classes.

Now consider the functor \ $T=(T_{i})_{i=0}^{\infty}:\CA{R}\to\gr\RM$, \ 
defined on cofree coalgebras \ by \ $T(\Gs)=\Vs$ \ if \ $\Gs=G(\Vs)$ \ (we can 
extend this by $0$-th derived functors to all of \ $\CA{R}$, \ if we wish). \ 
If \ $\Aus=\HRs{\Qu}$, \ then \ $T\Aus=\pi_{\ast}\Qu$, \ so \ 
$\pi^{k}\pi_{\ast}\Qu$ \ is just the $k$-th derived functor of $T$ applied to \ 
$\Ks$, \ denoted by \ $(L^{k}T)\Ks\in\gr\RM$ \ (cf.\ \cite[\S 7.8]{BlaN}), \ 
and it makes sense to say that \ $\Ks$ \ has \emph{projective dimension} 
$\leq n$ \ if \ $(L^{k}T)\Ks=0$ \ for \ $k>n$. \ For each \ $i\geq 0$, \ the 
functor \ $T_{i}$ \ has \emph{degree} $i$, in the sense dual to 
\cite[2.3.2]{BlaD}, so by (the dual of) \cite[Thm.\ 3.1]{BlaD} we have \ 
$\pi^{k}\pi_{s}\Qu=0$ \ for \ $k\geq s$. \ Then another convergence result of 
Bousfield's, namely, \cite[Thm.\ 3.4]{BousHS}, yields:

%
%
\begin{prop}\label{pfour}\stepcounter{subsection}
For \ $R=\Fp$, \ an unstable coalgebra \ $\Ks\in\CA{R}$ \ of
finite projective dimension is realizable as the homology of some space \ 
$\X\in\Ta$ \ if and only if \ $\Ks$ \ has a coherently vanishing sequence 
of characteristic classes.
\end{prop}

However, it is not clear on the face of it whether unstable coalgebras can 
\textit{ever} have non-trivial finite projective dimension (compare 
\cite[Thm.\ 4.3]{LMargH}).

\begin{remark}\label{ruar}\stepcounter{subsection}
As noted in the introduction, when \ $R=\Fp$, \ Theorem \ref{ttwo} (and perhaps
also Proposition \ref{pfour}) provides a way of constructing small, even
minimal, ``unstable Adams resolutions'' \ $\Qu$ \ of a given (simply-connected) 
space $\X$, which could be used in computing the Bousfield-Kan spectral 
sequence of \cite{BKaS} for \ $\pi_{\ast}(R_{\infty}\X)$. \ In particular, 
when $\X$ is of finite type, one can choose \ $\Qu$ \ so that each space \ 
$\bQ^{n}$ \ is a finite-type product of copies of \ $\EM{\Fp}{k}$ \ 
(for various $k$).
\end{remark}
%
%
\sect{Distinguishing between realizations}
\label{cdr}

Another interesting question is how one can distinguish between 
non-homotopy equivalent realizations \ $\X,\Y\in \Sa$ \ of a given unstable
coalgebra \ $\Ks$; \ we shall try to do this in terms of different realizations \ 
$\Qu$, \ $\Tu\in c\Sa$ \ of a fixed cosimplicial cofree resolution \ 
$\Ks\to\Aus$, \ where we assume to begin with that \ $\Ks$ \ is in fact 
realizable. Our goal is to find \emph{necessary} conditions in order for two 
realizations \ $\Qu$ \ and \ $\Tu$ \ to yield homotopy equivalent total spaces 
(compare \cite[Thm.\ 4.21]{BlaAI}).

Again the key lies in the extension of coalgebras \eqref{esix}. Of 
course, we may assume that the characteristic class \ 
$\chi_{n}\in H^{n+2}(\Ks;\Sigma^{n}\Ks)$ \ vanishes, so that it is possible
to find various splittings of \eqref{esix} as a ``semi-direct product'', 
given by different choices of the lifting $\lambda$  in Figure \ref{fig3}. 
The difference between two such semi-direct products is represented by a
suitable cohomology class (compare \cite[\S 6]{KLeeC} and 
\cite[IV, \S 2]{MacH}), constructed as follows:

\begin{defn}\label{ddoc}\stepcounter{subsection}
Given two liftings \ $\lambda,\lambda':\HRs{B^{n}\Qu}\to\bAs{n+1}$ \ 
in Figure \ref{fig3} above \ -- \ determining extensions of \ $\tr{n}\Qu$ \ to \ 
$\tr{n+1}\Qu$) \ -- \ as in the proof of Theorem \ref{tone}, we may assume that 
the corresponding maps \ $\mu,\mu':\Sigma^{n}\Ks\to B^{n+1}\Aus$ \ vanish. We 
extend \ $\lambda$, \ $\lambda'$ \ as in \S \ref{scc} to coface maps \ 
$d^{0},d_{0}':\bQ^{n}\to\bQ^{n+1}$, \ define \ 
$\eta:\HRs{B^{n}\Qu}\to C^{n+1}\Aus$ \ with \ $\eta\circ\bar{q}^{n}=0$, \ and 
extend to a map of \ $\As^{0}$-algebras \ $\zeta:\Sigma^{n}\Ks\to C^{n+1}\Aus$ \ 
(again, as in the proof of Theorem \ref{tone}). \ Again \ 
$\qu{n+1}{A^{n+1}}\circ \zeta=0$, \ so \ $\zeta$ is a cocycle 
in \ $\Hom_{\CA{\As^{0}}}(\Sigma\Ks,C^{\ast}\Aus)$, \ representing a cohomology 
class \ $\delta_{\lambda,\lambda'}\in H^{n+1}(\Ks,\Sigma^{n}\Js)$, \ which we 
call \ the \emph{difference obstruction} for the corresponding Postnikov 
sections \ $\bP{n}\Qu[\lambda]$ \ and \ $\bP{n}\Qu[\lambda']$ \ (in the notation 
of \S \ref{apost}).
\end{defn}

\begin{remark}\label{rdo}\stepcounter{subsection}
Again, by standard arguments this cohomology class is independent of the 
specific algebraic resolution \ $\Ks\to\Aus$ \ in \ $c\CA{R}$. \ 
Now assume that \ $\X,\Y\in\Sa$ \ are two (different) realizations of \ $\Ks$, \ 
with \ $\Qu$ \ and \ $\Tu$ \ respectively cosimplicial spaces realizing \ $\Aus$ \ 
(so that \ $\Qu\simeq \co{\X}$ \ and \ $\Tu\simeq \co{Y}$ \ in the resolution 
model category \ $c\Sa$ \ of \S \ref{srmc}), with the same \ $(n-1)$-type (that 
is, \ $\bP{n-1}\Qu\simeq\bP{n-1}\Tu$, \ so in particular we can assume that \ 
$\tr{n}\Qu=\tr{n}\Tu$). \ Then we can choose \ $\lambda$ \ and \ $\lambda'$ \ 
so that \ $\bP{n}\Qu[\lambda]=\bP{n}\Qu$, \ 
$\bP{n}\Qu[\lambda']=\bP{n}\Tu$ \ -- \ and thus \ $\delta_{\lambda,\lambda'}$ \ 
depends only on the $n$-type of \ $\Qu$ \ and \ $\Tu$, \ respectively, so in 
particular only on the homotopy types of $\X$ and $\Y$.
\end{remark}
%
%
\begin{thm}\label{tthree}\stepcounter{subsection}
If \ $\delta_{n}=0$ \ in \ $H^{n+1}(\Ks;\Sigma^{n}\Ks)$, \ 
then \ $\bP{n}\Qu[\lambda]\simeq\bP{n}\Qu[\lambda']$ \ in the resolution model 
category structure.
\end{thm}

\begin{proof}
If \ $\delta_{n}=0$, \ there is a map of \ $\As^{0}$-comodules \ 
$\vartheta:\Sigma^{n}\Ks\to C^{n}\Aus$ \ such that \ 
$\dd{A}\circ q^{n}_{A}\circ\vartheta=\zeta$, \ and by the 
discussion in \S \ref{rqc} \ $\vartheta$ \ can be lifted to a map \ 
$\theta:\Sigma^{n}\Ks\to \As^{n}$ \ (actually factoring through \ 
$(\As^{n})_{\ca}\hra\As^{n}$. \ If we define a map of \ $\As^{0}$-comodules: \ 
$\varphi\DEF \qu{n}{Q}\circ\vartheta\circ\psi^{n}$, \ then \  
$\lambda'\circ\varphi:\HRs{B^{n}\Qu}\to\bAs{n+1}$ \ is just \
\setcounter{equation}{\value{thm}}\stepcounter{subsection}
\begin{equation}\label{eeleven}
\bar{d}^{0}_{A}\circ\varphi=\proj_{\bar{A}}\circ\eta=\lambda-\lambda'
\end{equation}
\setcounter{thm}{\value{equation}}

\noindent in the following diagram:

%
%
\begin{figure}[htbp]
\begin{picture}(350,160)(-5,5)
%
%
\put(205,120){$\HRs{B^{n}\Qu}$}
\put(275,125){\vector(1,0){25}}
\put(300,125){\vector(1,0){3}}
\put(282,130){$\psi^{n}$}
\put(307,120){$\Coker(\bar{q}^{n})\cong\Sigma^{n}\Ks$}
%

%
\put(335,135){\line(-1,1){20}}
\put(315,155){\line(-1,0){120}}
\put(195,155){\vector(-2,-1){180}}
\put(105,117){$\theta$}
\put(197,154){\vector(-3,-4){68}}
\put(147,100){$\vartheta$}
\put(240,112){\vector(0,-1){45}}
\put(243,90){$\varphi$}
\put(260,114){\vector(4,-3){68}}
\put(295,95){$\eta$}
\put(340,112){\vector(0,-1){45}}
\put(343,90){$\zeta$}
%
%
\put(0,50){$\HRs{\bQ^{n}}$}
\put(57,55){\vector(1,0){35}}
\put(62,62){$(p^{n}_{Q})_{\#}$}
\put(97,50){$\HRs{C^{n}\Qu}$}
\put(170,55){\vector(1,0){32}}
\put(173,62){$\qu{n}{Q}$}
\put(205,50){$\HRs{B^{n}\Qu}$}
\put(276,55){\vector(1,0){30}}
\put(282,62){$\dd{Q^{n}}$}
\put(308,50){$\HRs{C^{n+1}\Qu}\cong C^{n+1}\Aus$}
%
%
\put(270,45){\vector(4,-3){54}}
\put(274,23){$\bar{q}^{n}$}
\put(360,42){\vector(0,-1){30}}
\put(314,30){$(\proj_{\bar{Q}})_{\#}$}
\put(420,42){\vector(0,-1){30}}
\put(389,30){$\proj_{\bar{A}}$}
\put(325,0){$\HRs{\bar{\bQ}^{n+1}}\cong\bAs{n+1}$}
\end{picture}
\caption[fig1]{}
\label{fig5}
\end{figure}

Note that because \ $\As^{n}$ \ is cofree, \ we can realize \ 
$\theta\circ\psi^{n}$ \ by a map  \ $f:B^{n}\Qu\to\bQ^{n}$ \ in \ $\Sa$, \ so \ 
$\varphi$ \ is realized by \ 
$q^{n}_{Q}\circ p^{n}_{Q}\circ f:B^{n}\Qu\to B^{n}\Qu$.

We may take the simplicial GEM \ $\bQ^{n}$ \ to be a simplicial 
$R$-module, with \ $\nu:\bQ^{n}\times\bQ^{n}\to\bQ^{n}$ \ the addition map, \ and
define \ $g:\bQ^{n}\to\bQ^{n}$ \ to be the composite \ 
$\nu \circ(\Id\top (f\circ p^{n}\circ q^{n}))$. \ For every \ $0\leq i\leq n$ \ 
we have \ $g\circ d^{i}=d^{i}:\bQ^{n-1}\to\bQ^{n}$, \ so \ 
$p^{n}\circ q^{n}\circ g$ \ induces a map \ $h:B^{n}\Qu\to B^{n}\Qu$, \ with \ 
\setcounter{equation}{\value{thm}}\stepcounter{subsection}
\begin{equation}\label{etwelve}
h_{\#}(\alpha)=\alpha+\varphi(\alpha)\hsp \text{for} \ \alpha\in\HRs{B^{n}\Qu},  
\end{equation}
\setcounter{thm}{\value{equation}}

\noindent and thus by \eqref{eeleven} \ the following diagram in \ $\CA{R}$ \ 
commutes:

%
%
\begin{picture}(360,100)(-20,-16)
%
%
\put(0,50){$B^{n}\Aus=B^{n}\HRs{\Qu}$}
\put(120,55){\vector(1,0){35}}
\put(130,60){$\bar{q}^{n}$}
\put(160,50){$\HRs{B^{n}\Qu}$}
\put(235,55){\vector(1,0){35}}
\put(245,60){$\lambda$}
\put(275,50){$\HRs{\bar{\bQ}^{n+1}}=\bAs{n+1}$}
%
%
\put(75,44){\vector(0,-1){27}}
\put(79,28){$id$}
\put(190,45){\vector(0,-1){29}}
\put(193,28){$h_{\#}$}
\put(310,43){\vector(0,-1){27}}
\put(315,28){$id$}
%
%
\put(0,0){$B^{n}\Aus=B^{n}\HRs{\Qu}$}
\put(120,5){\vector(1,0){35}}
\put(130,10){$\bar{q}^{n}$}
\put(160,0){$\HRs{B^{n}\Qu}$}
\put(235,5){\vector(1,0){35}}
\put(245,10){$\lambda'$}
\put(275,0){$\HRs{\bar{\bQ}^{n+1}}=\bAs{n+1}$}
\end{picture}

\noindent which yields a map of \ $(n+1)$-truncated objects \ 
$\tr{n+1}\Qu[\lambda]\to \tr{n+1}\Qu[\lambda']$, \ or equivalently, a map \ 
$\rho:\bP{n}\Qu[\lambda]\to\bP{n}\Qu[\lambda']$, \ which clearly induces an
isomorphism in \ $\pi^{k}\HRs{-}$ \ for \ $k\leq n+1$.

Note that for any choice of \ $\lambda$ \ we have \ 
$\pi^{n+2}\bP{n}\Qu[\lambda]\cong\Image(\partial_{n+2})\cong\Sigma^{n+1}\Ks$, \ 
by \eqref{enine}. Since \ $\psi^{n}\circ\qu{n}{Q}=0$ \ by \eqref{esix} and Figure 
\ref{fig2}, \ we have \ $\psi^{n}\circ\varphi=0$, \ so by \eqref{etwelve} \ 
$\psi^{n}\circ h_{\#}=\psi^{n}$, \ and since the following square commutes:

%
%
\begin{picture}(100,90)(-110,-15)
%
%
\put(0,50){$\HRs{B^{n+1}\Qu[\lambda]}$}
\put(95,55){\vector(1,0){45}}
\put(108,60){$\partial_{n+1}$}
\put(145,50){$\HRs{B^{n}\Qu}$}
%
%
\put(40,45){\vector(0,-1){30}}
\put(23,27){$\rho_{\#}$}
\put(180,43){\vector(0,-1){30}}
\put(183,27){$\rho_{\#}=h_{\#}$}
%
%
\put(0,0){$\HRs{B^{n+1}\Qu[\lambda']}$}
\put(100,5){\vector(1,0){40}}
\put(110,10){$\partial_{n+1}$}
\put(145,0){$\HRs{B^{n}\Qu}$}
\end{picture}

\noindent and \ $h_{\#}\rest{\Image(\partial_{n+1})}$ \ is an isomorphism, so is
$h_{\#}\rest{\Image(\partial_{n+2})}$, \ by \eqref{eeight}. Thus \ 
$\pi^{n+2}(\rho_{\#})$ \ is an isomorphism, too, so $\rho$ is a weak equivalence
by \eqref{enine} for \ $\bP{n}\Qu$.
\end{proof}

%
%
\begin{thm}\label{tfour}\stepcounter{subsection}
For \ $R=\Fp$ \ or \ $\Q$, \ assume \ $\X,\Y\in\Sa$ \ are two $R$-good 
realizations of a given unstable coalgebra \ $\Ks$ of finite type, which is 
either simply-connected or of finite projective dimension. \ 
If the difference obstructions for $\X$ and $\Y$ (\S \ref{rdo}) all vanish,
then $\X$ and $\Y$ are $R$-equivalent (i.e., $R_{\infty}\X\simeq R_{\infty}\Y$).
\end{thm}

\begin{proof}
When \ $\Ks$ \ is of finite type, we can choose a cosimplicial resolution \ 
$\Aus\in c\CA{R}$, \ with a CW basis in which each \ $\bAs{n}$ \ (and thus each \ 
$\As^{n}$) \ is of finite type. Let \ $\Qu$ \ and \ $\Tu$ \ be cosimplicial 
spaces realizing \ $\Aus$, \ which are resolutions  (in the sense of \S 
\ref{srmc}) of $\X$ and $\Y$ respectively, as in \S \ref{rdo}. By the Theorem 
\ref{ttwo} (resp.,\  Proposition \ref{pfour}), we have \ 
$\HRs{\Tot\Qu}\cong\Ks\cong\HRs{\Tot\Tu}$. \ 

Let \ $R^{\Delta}\Wu\in c\Sa$ \ denote the diagonal of the bicosimplicial 
space \ $R\W^{\bullet\bullet}$ \ obtained from a given cosimplicial space \ 
$\Wu\in c\Sa$ \ by applying the Bousfield-Kan $R$-resolution functor
(\cite[I, \S 4.1]{BKaH}) dimensionwise to \ $\Wu$. \ 
By Theorem \ref{tthree} there is a map of cosimplicial spaces \ 
$\rho:\Qu\to\Tu$ \ which is a weak equivalence in \ $c\Sa$, \ so induces an
isomorphism in the \ $E^{2}$-terms of the homology spectral sequences for \ 
$\Qu$ \ and  \ $\Tu$. \ Since \ $\HRs{\bQ^{n}}$ \ and \ $\HRs{\bT^{n}}$ \ are of 
finite type for each \ $n\geq 0$, \ by \cite[Thm.\ 9.1]{ShipC}, \ $\rho$ induces a
homotopy equivalence \ $\Tot R^{\Delta}\Qu\xra{\simeq}\Tot R^{\Delta}\Tu$ \ 
(and similarly \ $\Tot\Qu\to R_{\infty}\X$ \ and \ $\Tot\Tu\to R_{\infty}\Y$). \ 

However, for each \ $n\geq 0$, \ the spaces \ $\bQ^{n}$ \ and \ $\bT^{n}$ \ are 
$R$-GEMs, so they are $R$-complete (cf.\ \cite[V, 3.3]{BKaH}), and thus \ 
$\Tot R^{\Delta}\Qu\simeq\Tot\Qu$ \ by \cite[Thm.\ 10.2]{ShipC}, and similarly 
for \ $\Tu$, \ so we find that $\X$ and $\Y$ are indeed $R$-equivalent.
\end{proof}

\begin{remark}\label{rrat}\stepcounter{subsection}
Shipley's theorems, in \cite{ShipC}, were originally stated for \ $R=\Fp$; \ 
when \ $R=\Q$ \ it is no longer true that all relevant homotopy groups are 
finite.  However, they are finite dimensional vector spaces over $\Q$, \ so 
\cite[IX, \S 3]{BKaH}, and the rest of Shipley's arguments, still apply.

As noted in the proof of Theorem \ref{ttwo}, when \ $R=\Q$ \ the only problem
of interest is to distinguish between different realizations of a given \ 
$R$-(co)algebra; the  obstruction theory described here is just the vector-space 
dual of the theory for graded algebras over $\Q$ defined by Halperin and 
Stasheff in \cite{HStaO} (see also \cite{FelDT}).
\end{remark}

\begin{remark}\label{rdpt}\stepcounter{subsection}
When \ $R=\Fp$, \ Theorem \ref{tfour} can be thought of as providing a collection 
of algebraic invariants \ -- \ starting with the homology coalgebra \ 
$\His{\X}{\Fp}$  \ -- \ for distinguishing between $p$-types of spaces. 
As with the ordinary Postnikov systems and their $k$-invariants, these are not 
actually invariant, in the sense that distinct values (i.e., non-vanishing 
difference obstructions) do not guarantee distinct $p$-types. 

This approach is the Hilton-Eckmann dual of the theory described in \cite{BlaAI} 
or \cite{BGoeC} for distinguishing (integral) homotopy types, starting with the 
homotopy $\Pi$-algebra \ $\pi_{\ast}\X$, \ in terms of an analogous collection of 
cohomology classes. It is reasonable to expect a more general version of 
Theorem \ref{tfour} to hold, without the assumption of finite type, and for any \ 
$R\subseteq\Q$; \ but this would require a stronger convergence result than that
provided by \cite[\S 9-10]{ShipC}.

Perhaps it should be observed that many non-realization results proven in the past
(see Introduction) have used higher order cohomology operations; these are 
implicit in the Quillen cohomology cohomology classes of Theorems \ref{tone} and 
\ref{tthree}, and were made explicit in the $\Pi$-algebra analogue in 
\cite{BlaHH}. We hope to return to this point in the future.
\end{remark}


\begin{thebibliography}{ABC2}
%
\bibitem[A]{AdHI}
J.F.~Adams ,
``On the non-existence of elements of {Hopf} invariant one'',\hsm
\textit{Ann. Math. (2)} \textbf{72} (1960), No.\ 1, pp.\ 20-104.
%
\bibitem[AW]{AWiF}
J.F.~Adams \& C.W.~Wilkerson,
``Finite {$H$}-spaces and algebras over the {Steenrod} algebra'',\hsm
\textit{Ann. Math.} \textbf{111} (1980), pp.\ 95-143.
%
\bibitem[Ad]{AdemI}
J.~Adem,
``The iteration of the {Steenrod} squares in algebraic topology'',\hsm
\textit{Proc. Nat. Acad. Sci. USA} \textbf{38} (1952), pp.\ 720-726.
%
\bibitem[Ag]{AguRC}
J.~Aguad{\'{e}},
``Realizability of cohomology algebras: a survey'',\hsm
\textit{Pub. Mat. Univ. Aut. Barcelona} \textbf{26} (1982), No.\ 2, pp.\ 25-68.
%
\bibitem[ABN]{ABNotH1}
J.~Aguad{\'{e}}, C.~Broto, \& D.~Notbohm.
``Homotopy classification of some spaces with interesting cohomology and a 
conjecture of Cooke, Part I'',\hsm
\textit{Topology} \textbf{33} (1994), No.\ 3, pp.\ 455-492.
%
\bibitem[An]{AndG}
D.W.~Anderson,
``A generalization of the {Eilenberg}-{Moore} spectral sequence'',\hsm
\textit{Bull. AMS} \textbf{78} (1972), No.\ 5, pp.\ 784-786.
%
\bibitem[And]{AndrM}
M.~Andr\'{e},
\textit{M\'{e}thode Simpliciale en Alg\`{e}bre Homologique et Alg\`{e}bre 
Commutative},\hsm
Springer-\-Verlag \textit{Lec.\ Notes Math.} \textbf{32}, Berlin-\-New York, 1967.
%
\bibitem[BM]{BMadR}
N.A.~Baas \& I.H.~Madsen,
``On the realization of certain modules over the {Steenrod} algebra'',\hsm
\textit{Math.\ Scand.} \textbf{31} (1972), pp.\ 220-224.
%
\bibitem[Bl1]{BlaD}
D.~Blanc,
``Derived functors of graded algebras'',\hsm
\textit{J.\ Pure Appl.\ Alg.} \textbf{64} (1990) No.\ 3, pp.\ 239-262.
%
\bibitem[Bl2]{BlaN}
D.~Blanc,
``New model categories from old'',\hsm
\textit{J.\ Pure \& Appl.\ Alg.} \textbf{109} (1996) No.\ 1, pp.\ 37-60.
%
\bibitem[Bl3]{BlaHH}
D.\ Blanc,
``Higher homotopy operations and the realizability of homotopy groups'',\hsm 
\textit{Proc.\ Lond.\ Math.\ Soc.\ (3)} \textbf{70} (1995), pp.\ 214-240.
%
\bibitem[Bl4]{BlaC}
D.~Blanc,
``CW simplicial resolutions of spaces, with an application to loop spaces'',\hsm 
to appear in \textit{Topology \& Appl.}.
%
\bibitem[Bl5]{BlaAI}
D.~Blanc,
``Algebraic invariants for homotopy types'',\hsm 
to appear in \textit{Math.\ Proc.\ Camb.\ Philos.\ Soc.}.
%
\bibitem[BDG]{BDGoeC}
D.~Blanc, W.G.~Dwyer, \& P.G.~Goerss,
``Cohomology invariants for simplicial spaces'',\hsm
preprint 1999.
%
\bibitem[BS]{BStG}
D.~Blanc \& C.R.~Stover,
``A generalized {Grothendieck} spectral sequence'',\hsm
in N.\ Ray \& G.\ Walker, eds., \textit{Adams Memorial Symposium on
  Algebraic Topology, Vol.\ 1}, Lond.\ Math.\ Soc.\ Lec.\
Notes Ser.\ \textbf{175}, Cambridge U. Press, Cambridge, 1992,  pp.\ 145-161. 
%
\bibitem[Bo1]{BousHS}
A.K.~Bousfield,
``On the homology spectral sequence of a cosimplicial space'',\hsm
\textit{Amer.\ J.\ Math.} \textbf{109} (1987), No.\ 2, pp.\ 361-394.
%
\bibitem[Bo2]{BousH}
A.K.~Bousfield,
``Homotopy spectral sequences and obstructions'',\hsm
\textit{Isr.\ J.\ Math.} \textbf{66} (1989), Nos.\ 1-3, pp.\ 54-104.
%
\bibitem[BC]{BCurS}
A.K.~Bousfield \& E.B.~Curtis,
``A spectral sequence for the homotopy of nice spaces'',\hsm
\textit{Trans.\ AMS} \textbf{151} (1970), pp.\ 457-478.
%
\bibitem[BK1]{BKaH}
A.K.~Bousfield \& D.M.~Kan,
\textit{Homotopy Limits, Completions, and Localizations},\hsm
Springer-\-Verlag \textit{Lec.\ Notes Math.} \textbf{304}, 
Berlin-\-New York, 1972.
%
\bibitem[BK2]{BKaS}
A.K.~Bousfield \& D.M.~Kan,
``The homotopy spectral sequence of a space with coefficients in a ring'',\hsm
\textit{Topology} \textbf{11} (1972), pp.\ 79-106.
%
\bibitem[BG]{BGitS}
E.H.\ Brown, Jr.\ \& S.~Gitler,
``A spectrum whose cohomology is a certain cyclic module over the Steenrod 
algebra'',\hsm
\textit{Topology} \textbf{12} (1973) No.\ 3, pp.\ 283-293.
%
\bibitem[BP]{BPetS}
E.H.\ Brown, Jr.\ \& F.P.~Peterson,
``A spectrum whose $Z_{p}$ cohomology is the algebra of reduced $p$-th powers,\hsm
\textit{Topology} \textbf{5} (1966) No.\ 2, pp.\ 149-154.
%
\bibitem[CE]{CEwiR}
A.~Clark \& J.H.~Ewing,
``The realization of polynomial algebras as cohomology rings'',\hsm
\textit{Pac. J. Math.} \textbf{50} (1974), No.\ 2, pp.\ 425-434.
%
\bibitem[CS]{CSmiRM}
G.E.~Cooke \& L.E.~Smith,
``On realizing modules over the {Steenrod} algebra'',\hsm
\textit{J. Pure Appl. Alg.} \textbf{13} (1978), No.\ 1, pp.\ 71-100.
%
\bibitem[Do]{DoH}
A.~Dold,
``Homology of symmetric products and other functors of complexes'',\hsm 
\textit{Ann.\ Math.\ (2)} \textbf{68} (1958), pp.\ 54-80.
%
\bibitem[DKW]{DKWinsC}
J.~Duflot, N.J.~Kuhn, \& M.W.~Winstead,
``A classification of polynomial algebras over the {Steenrod} algebra'',\hsm
\textit{Comm. Math. Helv.} \textbf{68} (1993), No.\ 4, pp.\ 622-632.
%
\bibitem[DKS1]{DKStE} 
W.G.~Dwyer, D.M.~Kan, \& C.R.~Stover, 
``An $E^{2}$ model category structure for pointed simplicial spaces'',\hsm 
\textit{J.\ Pure \& Appl.\ Alg.} \textbf{90} (1993) No.\ 2, pp.\ 137-152.
%
\bibitem[DKS2]{DKStB} 
W.G.~Dwyer, D.M.~Kan, \& C.R.~Stover, 
``The bigraded homotopy groups $\pi_{i,j}X$ of a pointed simplicial space'',\hsm 
\textit{J.\ Pure Appl.\ Alg.} \textbf{103} (1995), No.\ 2, pp.\ 167-188.
%
\bibitem[DMW]{DMWilH}
W.G. Dwyer, H.R.~Miller, \& C.C.~Wilkerson,
``Homotopical uniqueness of classifying spaces'',\hsm
\textit{Topology} \textbf{31} (1992), No.\ 1, pp.\ 29-45.
%
\bibitem[DW]{DWilF}
W.G. Dwyer \& C.C.~Wilkerson,
``A new finite loop space at the prime two'',\hsm
\textit{J. AMS} \textbf{6} (1993) No.\ 1, pp.\ 37-64.
%
\bibitem[EM]{EMoHF}
S.~Eilenberg \& J.C.~Moore,
``Homology and fibrations, I: Coalgebras, cotensor product and its derived 
functors'',\hsm
\textit{Comm.\ Math.\ Helv.} \textbf{40} (1966), pp.\ 199-236.
%
\bibitem[F]{FelDT}
Y.~F\'{e}lix,
\textit{D\`{e}nombrement des types de $k$-homotopie: th\'{e}orie de 
la d\'{e}formation},\hsm
Mem.\ Soc.\ Math.\ France \textbf{3}, Paris, 1980.
%
\bibitem[G]{GoeHH}
P.G.~Goerss,
``The homology of homotopy inverse limits'',\hsm
\textit{J.\ Pure Appl.\ Alg.} \textbf{111} (1996), No.\ 1, pp.\ 83-122.
%
\bibitem[GH]{GHopR}
P.G.~Goerss \& M.J.~Hopkins,
``Resolutions in model categories'',\hsm
preprint 1998.
%
\bibitem[HS]{HStaO}
S.~Halperin \& J.D.~Stasheff,
``Obstructions to homotopy equivalences'',\hsm
\textit{Adv.\ in Math.} \textbf{32} (1979) No.\ 3, pp.\ 233-279.
%
\bibitem[Ha]{HarpC1}
J.R.~Harper,
``On the construction of mod $p$ $H$-spaces'',\hsm
in R.J.~Milgram, ed., \textit{Algebraic and Geometric Topology, Part 2 
(Stanford,CA, 1976)}, Proc.\ Symp.\ Pure Math. \textbf{32}, AMS,  
Providence, RI, 1978, pp.\ 207-214.
%
\bibitem[Ho]{HopfT}
H.~Hopf,
``\"{U}ber die Topologie der Gruppen-Mannigfaltkeiten und ihre 
Verallgemeinerungen'',\hsm
\textit{Ann.\ Math.\ (2)} \textbf{42} (1941), pp.\ 22-52.
%
\bibitem[Ka]{KanR}
D.M.~Kan,
``A relation between $CW$-complexes and free c.s.s.\ groups'',\hsm 
\textit{Am.\ J.\ Math.} \textbf{81} (1959), pp.\ 512-528.
%
\bibitem[Ku]{KuhnTR}
N.J.~Kuhn,
``On topologically realizing modules over the {Steenrod} algebra'',\hsm
\textit{Ann. Math., Ser. 2} \textbf{141} (1995), No.\ 2, pp.\ 321-347.
%
\bibitem[L]{KLeeC}
K.~Lee",
``Cosimplicial cohomology of coalgebras'',\hsm
\textit{Nagoya Math.\ J.} \textbf{47} (1972), pp.\ 199-226.
%
\bibitem[LM]{LMargH}
T.Y.~Lin \& H.R.~Margolis,
\newblock 
``Homological aspects of modules over the {Steenrod} algebra'',\hsm
\textit{J.\ Pure Appl.\ Alg.} \textbf{9} (1976/77), No.\ 2, pp.\ 121-129.
%
\bibitem[M1]{MacH} 
S.~Mac~Lane, 
\textit{Homology},\hsm 
Springer-Verlag \textit{Grund.\ math.\ Wissens.} \textbf{114}, 
Berlin-\-New York  1963.
%
\bibitem[M2]{MacC}
S.~Mac Lane,
\textit{Categories for the Working Mathematician},\hsm
Springer-\-Verlag \textit{Grad.\ Texts in Math.} \textbf{5}, 
Berlin-\-New York, 1971.
%
\bibitem[M]{MayS}
J.P.~May,
\textit{Simplicial Objects in Algebraic Topology},\hsm
U.\ Chicago Press, Chicago-\-London, 1967.
%
\bibitem[Ma]{MargS}
H.R.~Margolis,
\textit{Spectra and the Steenrod Algebra: \ Modules over the Steenrod Algebra 
and the Stable Homotopy Category},\hsm
North-Holland, Amsterdam-\-New York, 1983.
%
\bibitem[Ml]{MilSC} 
H.~Miller, 
``Correction to `The Sullivan conjecture on maps from classifying spaces' ''\hsm 
\textit{Ann.\ of Math.} \textbf{121} (1985), pp.\ 605-609.
%
\bibitem[Mn]{MilnS}
J.W.~Milnor,
``The Steenrod algebra and its dual'',\hsm
\textit{Ann.\ Math.\ (2)} \textbf{67} (1958), pp.\ 150-171.
%
\bibitem[MM]{MMoorH}
J.W.~Milnor \& J.C.~Moore,
``On the structure of Hopf algebras'',\hsm
\textit{Ann.\ Math.\ (2)} \textbf{81} (1965), pp.\ 211-264.
%
\bibitem[Q1]{QuH}
D.G.~Quillen,
\textit{Homotopical Algebra},\hsm
Springer-\-Verlag \textit{Lec.\ Notes Math.} \textbf{20}, 
Berlin-\-New York, 1963.
%
\bibitem[Q2]{QuR}
D.G.~Quillen,
``Rational homotopy theory'',\hsm
\textit{Ann.\ Math.\/} \textbf{90} (1969) No.\ 2, pp.\ 205-295.
%
\bibitem[Q3]{QuC} 
D.G.~Quillen, 
``On the (co-)homology of commutative rings'',\hsm 
\textit{Applications of Categorical Algebra}, \ Proc.\ Symp.\ Pure Math.\ 
\textbf{17}, AMS, Providence, RI, 1970, pp.\ 65-87.
%
\bibitem[R]{RecS}
D.L.~Rector, 
``Steenrod operations in the Eilenberg-Moore spectral sequence'',\hsm
\textit{Comm. Math. Helv.} \textbf{45} (1970), pp.\ 540-552.
%
\bibitem[Sc1]{SchwU}
L.~Schwartz,
\textit{Unstable Modules over the Steenrod Algebra and Sullivan's Fixed Point 
Set Conjecture},\hsm
U. Chicago Press, Chicago-\-London, 1994.
%
\bibitem[Sc2]{SchwN}
L.~Schwartz,
``{`{A}} propos de la conjecture de non-r\'{e}alisation due \`{a} N. Kuhn'',\hsm
\textit{Inv.\ Math.} \textbf{134} (1998), No.\ 1, pp.\ 211-227.
%
\bibitem[Sh]{ShipC}
B.E.~Shipley, 
``Convergence of the homology spectral sequence of a cosimplicial space'',\hsm
\textit{Amer. J. Math.} \textbf{118} (1996), No.\ 1, pp.\ 179-207.
%
\bibitem[Sm]{LSmiRS}
L.E.~Smith,
``On the realization and classification of symmetric algebras as cohomology 
rings'',\hsm
\textit{Proc.\ AMS} \textbf{87} (1983), No.\ 1, pp.\ 144-148.
%
\bibitem[SS]{SSwitR}
L.E.~Smith \& R.M.~Switzer,
``Realizability and nonrealizability of Dickson algebras as cohomology 
rings'',\hsm 
\textit{Proc. AMS} \textbf{89} (1983), No.\ 2, pp.\ 303-313.
%
\bibitem[St1]{SteCOP}
N.E.~Steenrod,
``Cohomology operations'',\hsm
in J.~Adem et~al., eds, 
\textit{Symposium internacional de topolog\'{\i}a algebraica}, 
UNAM/UNESCO, Mexico City, 1958, pp.\ 165-185. 
%
\bibitem[St2]{SteCA}
N.E.~Steenrod,
``The cohomology algebra of a space'',\hsm
\textit{Ens. Math.} \textbf{7} (1961), pp.\ 153-178.
%
\bibitem[ST]{STodS}
T.~Sugawara \& H.~Toda,
``Squaring operations in truncated polynomial algebras'',\hsm
\textit{Jap.\ J.\ Math.} \textbf{38} (1969), pp.\ 39-50.
%
\bibitem[Sw]{SweeH}
M.E.~Sweedler,
\textit{Hopf Algebras},\hsm W.A.~Benjamin, New York, 1969.
%
\bibitem[Ta]{TateH}
J.~Tate,
``Homology of noetherian rings and local rings'',\hsm
\textit{Ill.\ J.\ Math.} \textbf{1} (1957), pp.\ 14-27.
%
\bibitem[Th]{EThomS2}
E.~Thomas,
``Steenrod squares and $H$-spaces, II'',\hsm
\textit{Ann.\ Math.\ (2)} \textbf{81} (1965), pp.\ 473-495.
%
\bibitem[Wh]{GWhE}
G.W.~Whitehead,
\textit{Elements of Homotopy Theory},\hsm
Springer-\-Verlag \textit{Grad.\ Texts Math.} \textbf{61},
Berlin-\-New York, 1971.
%
\end{thebibliography}
\end{document}